\documentclass{article}

\title{Reflection Character Varieties in \(\mathrm{SL}_{\pm}(n,\mathbb R)\) for \(2\)-Perfect Labeled Truncation Polytopes}
\author{Suhyoung Choi and Yongho Seo\\
{\small{Department of Mathematical Sciences, KAIST}}\\
{\small{\href{mailto:schoi@math.kaist.ac.kr}{\texttt{schoi@math.kaist.ac.kr}}\quad
\href{mailto:seojm8@kaist.ac.kr}{\texttt{seojm8@kaist.ac.kr}}}}}
\date{}

\usepackage[T1]{fontenc}
\usepackage[utf8]{inputenc}
\usepackage[a4paper,margin=1.05in]{geometry}
\usepackage{graphicx}
\graphicspath{{figures/}}
\usepackage{amsthm,amssymb}
\usepackage{mathtools}
\usepackage{mathrsfs}
\usepackage{enumitem}
\usepackage[backend=biber,style=numeric,sorting=nty,maxcitenames=2,maxbibnames=99]{biblatex}
\usepackage{hyperref}
\hypersetup{
  hidelinks,
  pdftitle={Reflection Character Varieties in SL±(n, R) for 2-Perfect Labeled Truncation Polytopes},
  pdfauthor={Suhyoung Choi and Yongho Seo},
  pdfsubject={Reflection character varieties of Coxeter groups associated with truncation d-polytopes},
  pdfkeywords={Coxeter groups, reflection character varieties, Cartan matrices, truncation polytopes, projective geometry}
}

\newtheoremstyle{extrabolddef}
{\topsep}
{\topsep}
{\normalfont}
{}
{\bfseries\boldmath}
{.}
{ }
{\thmname{#1}\thmnumber{ #2}\thmnote{ (#3)}}

\newtheorem{theorem}{Theorem}[section]
\newtheorem{lemma}[theorem]{Lemma}
\newtheorem{proposition}[theorem]{Proposition}
\newtheorem{corollary}[theorem]{Corollary}

\theoremstyle{extrabolddef}
\newtheorem{definition}[theorem]{Definition}
\newtheorem{remark}[theorem]{Remark}

\theoremstyle{definition}
\newtheorem{example}[theorem]{Example}

\newcommand{\sslash}{\mathbin{/\mkern-6mu/}}
\DeclareMathOperator{\Hom}{Hom}
\DeclareMathOperator{\Span}{span}
\DeclareMathOperator{\diag}{diag}
\DeclareMathOperator{\rank}{rank}
\DeclareMathOperator{\Id}{Id}

\newcommand{\SLpm}[1]{\mathrm{SL}_{\pm}(#1,\mathbb{R})}
\newcommand{\PGL}[1]{\mathrm{PGL}(#1,\mathbb{R})}
\newcommand{\Homrefc}[1]{\operatorname{Hom}^{\mathrm{ref}}_{#1}}
\newcommand{\HomnVrc}[1]{\operatorname{Hom}^{\mathrm{ref}}_{#1,\mathrm{nVr}}}
\newcommand{\HomVrc}[1]{\operatorname{Hom}^{\mathrm{ref}}_{#1,\mathrm{Vr}}}
\newcommand{\Homstrictc}[1]{\operatorname{Hom}^{\mathrm{ref}}_{#1,\mathrm{str}}}
\newcommand{\Ccal}{\mathcal C}
\newcommand{\Homblockc}[1]{\operatorname{Hom}^{\mathrm{ref}}_{#1,\mathrm{blk}}}
\newcommand{\Rcal}{\mathscr R}
\newcommand{\Lcal}{\mathscr L}
\newcommand{\Ecal}{\mathscr E}
\newcommand{\Xn}{X_{n}}
\newcommand{\Xin}{\Xi_{n,\mathcal G}}
\newcommand{\Xref}[1]{X^{\mathrm{ref}}_{#1}}
\newcommand{\Xstr}[1]{X^{\mathrm{ref}}_{#1,\mathrm{str}}}
\newcommand{\XnnVr}{X^{\mathrm{ref}}_{n,\mathrm{nVr}}}
\newcommand{\XnVr}{X^{\mathrm{ref}}_{n,\mathrm{Vr}}}
\newcommand{\Xnblock}{X^{\mathrm{ref}}_{n,\mathrm{blk}}}
\newcommand{\Homranknc}[1]{\operatorname{Hom}^{\mathrm{ref},(n)}_{#1,\mathrm{str}}}

\begin{document}
	\maketitle

\begin{abstract}
	A truncation polytope is a combinatorial polytope represented by a convex polytope obtained from a simplex by a finite sequence of vertex truncations.
	Marquis studied the deformation spaces of labeled truncation \(3\)-polytopes, and Choi-Lee-Marquis studied the \(2\)-perfect case in higher dimensions.
	In this paper, we consider Coxeter groups associated with non-simplex \(2\)-perfect labeled truncation \(d\)-polytopes with connected Coxeter diagrams, where \(d\ge3\).
	For \(n\ge d+1\), we study their reflection-generated representations in \(\mathrm{SL}_{\pm}(n,\mathbb R)\) satisfying Vinberg's conditions.
	Their images in the real GIT quotients of the full representation spaces under conjugation define the reflection character varieties.
	These varieties and their subspaces of strict characters are semialgebraic, and we determine their common dimension whenever strict characters exist.
	We also decompose these subspaces by Cartan matrix rank, show that the nonempty pieces are smooth manifolds, and compute their dimensions explicitly.
\end{abstract}

\begingroup
\setlength{\baselineskip}{11pt}
\tableofcontents
\endgroup

\section{Introduction}\label{sec:intro}
A combinatorial \(d\)-polytope is a combinatorial equivalence class of convex \(d\)-polytopes.
A \emph{truncation \(d\)-polytope} is a combinatorial \(d\)-polytope represented by a convex polytope obtained from a \(d\)-simplex by a finite sequence of vertex truncations.
A \emph{facet} is a face of codimension \(1\), and a \emph{ridge} is a face of codimension \(2\).
Let \(\mathcal G\) be a labeled truncation polytope of dimension \(d\ge3\), with each ridge labeled by an integer at least \(2\).
Let \(S\) be the facet set of \(\mathcal G\), and set \(f=|S|\).
Let \(e\) be the number of ridges of \(\mathcal G\), and let \(e_2\) be the number of ridges labeled by \(2\).

Two distinct facets are called adjacent if their intersection is a ridge.
The associated Coxeter matrix \(M=(m_{st})_{s,t\in S}\) has diagonal entries \(m_{ss}=1\).
For distinct \(s,t\in S\), the entry \(m_{st}\) is the label of their common ridge when they are adjacent, and is \(\infty\) otherwise.

The associated Coxeter group is
\[
\pi_{\mathcal G}
=
\left\langle S \,\middle|\,
s^2=1\ \text{for }s\in S,\quad
(st)^{m_{st}}=1\ \text{for }s\ne t\text{ with }m_{st}<\infty
\right\rangle,
\]
where we use the same symbols for facets and the corresponding generators.
The Coxeter diagram of \(\pi_{\mathcal G}\) is the labeled graph with node set \(S\).
Two distinct nodes \(s,t\) are joined by an edge exactly when \(m_{st}\in\{3,4,\ldots,\infty\}\).
Each such edge has label \(m_{st}\), with the label \(3\) omitted.

The labeled polytope \(\mathcal G\) is \emph{\(2\)-perfect} if, for each edge \(E\) of \(\mathcal G\), the subgroup of \(\pi_{\mathcal G}\) generated by the facets containing \(E\) is finite~\cite[Section~1.2]{choi2022deformation}.
This condition is automatic for \(d=3\), since the two facets containing any edge generate a finite dihedral subgroup.
We assume that \(\mathcal G\) is \(2\)-perfect.

A properly convex \(d\)-polytope \(P\) in \(\mathbb S^d\) together with the projective reflections in the facets of \(P\) is a \emph{projective Coxeter polytope}, or simply \emph{Coxeter polytope}, provided that if \(\Gamma_P\) is the subgroup of \(\SLpm{d+1}\) generated by those reflections, then
\[
\operatorname{Int}(P)\cap\gamma\cdot\operatorname{Int}(P)=\varnothing
\quad\text{for each non-identity element }\gamma\in\Gamma_P,
\]
where \(\operatorname{Int}(P)\) denotes the interior of \(P\)~\cite[Section~1.2]{choi2022deformation}.
It follows from work of Vinberg~\cite[Theorem~2]{vinberg1971discrete} that the interior \(\Omega_P\) of the union of all \(\Gamma_P\)-translates of \(P\) is a convex domain of \(\mathbb S^d\) and that the group \(\Gamma_P\) acts properly discontinuously on \(\Omega_P\).
The quotient \(\widehat P=\Omega_P/\Gamma_P\) is then a convex real projective Coxeter orbifold.

A \emph{Coxeter realization} of \(\mathcal G\) in \(\mathbb S^d\) is a projective Coxeter \(d\)-polytope \(P\) whose face poset is identified with that of \(\mathcal G\) by an isomorphism preserving the facet indices and Coxeter labels (see Definition~\ref{def:coxeter-realization}).
We say that \(\mathcal G\) is \emph{Coxeter realizable} in \(\mathbb S^d\) if it admits such a realization.
For such a realization \(P\), the \(2\)-perfectness of \(\mathcal G\) implies that the underlying space of \(\widehat P\) is homeomorphic to \(P\) with a possibly empty set of vertices removed~\cite[Theorem~2.2\textup{(5)}]{choi2022deformation}.
The singular locus of \(\widehat P\) corresponds to the boundary of \(P\) with the same vertices removed.

For \(d=3\), Marquis~\cite[Theorem~3.16]{marquis2010espace} proved that the deformation space of Coxeter polytopes realizing \(\mathcal G\), if nonempty, is a finite disjoint union of open cells of the same dimension.
For \(d\ge4\), Choi-Lee-Marquis~\cite[Theorem~A]{choi2022deformation} proved the same conclusion when \(\mathcal G\) is \(2\)-perfect and \(\pi_{\mathcal G}\) is irreducible and large.

The facet reflections of a Coxeter realization of \(\mathcal G\) determine a representation of \(\pi_{\mathcal G}\) in \(\SLpm{d+1}\).
In this paper, we study representations
\[
\rho:\pi_{\mathcal G}\longrightarrow\SLpm{n},
\qquad n\ge d+1,
\]
that send each standard generator to a reflection.

Let \(\rho\in\Hom(\pi_{\mathcal G},\SLpm{n})\).
A \emph{reflection pair} for \(\rho\) is a family \((\alpha_s,v_s)_{s\in S}\), with \(\alpha_s\in(\mathbb R^n)^*\) and \(v_s\in\mathbb R^n\), such that
\[
\rho(s)=\Id-\alpha_s\otimes v_s,
\qquad
\alpha_s(v_s)=2
\quad\text{for every }s\in S.
\]
Here \((\alpha\otimes v)(x)=\alpha(x)v\).
The associated \emph{Cartan matrix} is \(C=(c_{st})_{s,t\in S}\), where \(c_{st}=\alpha_s(v_t)\).

The \emph{fundamental cone} associated with the reflection pair is
\[
\widetilde{\Delta}
=
\{x\in\mathbb R^n\mid \alpha_s(x)\le0\text{ for every }s\in S\}.
\]
Its image in the projective sphere \(\mathbb S^{n-1}\) is called the \emph{fundamental chamber}.

We call \(\rho\) \emph{\(\mathcal G\)-compatible} if it admits a reflection pair satisfying the following \emph{Vinberg conditions}:
\begin{enumerate}[label=(V\arabic*),labelindent=\parindent,leftmargin=*,itemsep=0.2ex,topsep=0.4ex]
	\item \(c_{st}=0\) for all \(s\ne t\) with \(m_{st}=2\).
	\item \(c_{st}<0\) for all \(s\ne t\) with \(m_{st}\ne2\).
	\item \(c_{st}c_{ts}=4\cos^2(\pi/m_{st})\) for all \(s\ne t\) with \(2<m_{st}<\infty\).
	\item \(c_{st}c_{ts}\ge4\) for all \(s\ne t\) with \(m_{st}=\infty\).
	\item \(\widetilde{\Delta}\) has nonempty interior in \(\mathbb R^n\).
\end{enumerate}
We denote the space of these representations by \(\Homrefc{\mathcal G}(\pi_{\mathcal G},\SLpm{n})\).
See Definition~\ref{def:vinberg-conditions} for the general formulation.

Let \(\rho\in\Homrefc{\mathcal G}(\pi_{\mathcal G},\SLpm{n})\), and choose a reflection pair \((\alpha_s,v_s)_{s\in S}\) for \(\rho\).
We call \(\rho\) \emph{Vinberg-reduced} if
\[
\bigcap_{s\in S}\ker\alpha_s=\{0\}.
\]
We call \(\rho\) \emph{strict} if one, and hence every, associated Cartan matrix \(C=(c_{st})_{s,t\in S}\) satisfies \(c_{st}c_{ts}>4\) whenever \(m_{st}=\infty\).
We denote the subspace of strict representations by \(\Homstrictc{\mathcal G}(\pi_{\mathcal G},\SLpm{n})\).

For \(1\le r\le n\), let \(\iota_r^n:\SLpm{r}\hookrightarrow\SLpm{n}\) be the standard block embedding defined by \(\iota_r^n(A)=\diag(A,\Id_{n-r})\).
If \(\rho_0:\pi_{\mathcal G}\to\SLpm{d+1}\) is the representation determined by a Coxeter realization of \(\mathcal G\), then its block extension \(\iota_{d+1}^n\circ\rho_0=\rho_0\oplus\Id_{\mathbb R^{n-d-1}}\) is \(\mathcal G\)-compatible and is reducible when \(n>d+1\).
We call a \(\mathcal G\)-compatible representation \(\rho:\pi_{\mathcal G}\to\SLpm{n}\) \emph{block-embedded} if it is \(\SLpm{n}\)-conjugate to \(\iota_{n-1}^n\circ\sigma\) for some \(\sigma\in\Homrefc{\mathcal G}(\pi_{\mathcal G},\SLpm{n-1})\).

The orbit space \(\Hom(\pi_{\mathcal G},\SLpm{n})/\SLpm{n}\) for the conjugation action need not be Hausdorff.
We instead use the real geometric invariant theory (GIT) quotient
\[
\begin{aligned}
	\Xn(\mathcal G)
	&:=\Hom(\pi_{\mathcal G},\SLpm{n})\sslash\SLpm{n},\\
	\Xin
	&:\Hom(\pi_{\mathcal G},\SLpm{n})\longrightarrow\Xn(\mathcal G),
\end{aligned}
\]
where \(\Xin\) is the quotient map.
For \(\rho\in\Hom(\pi_{\mathcal G},\SLpm{n})\), write \([\rho]:=\Xin(\rho)\).
Birkes~\cite[Theorem~5.2]{birkes1971orbits} and Richardson-Slodowy~\cite[Section~7.3]{richardson1990minimum} show that two representations have the same image under \(\Xin\) if and only if the closures of their \(\SLpm{n}\)-orbits meet.
Here orbit closures are taken in the usual topology.

We define the \emph{reflection character variety} and the \emph{strict reflection character variety} by
\[
\begin{aligned}
	\Xref{n}(\mathcal G)
	&:=\Xin\bigl(\Homrefc{\mathcal G}(\pi_{\mathcal G},\SLpm{n})\bigr),\\
	\Xstr{n}(\mathcal G)
	&:=\Xin\bigl(\Homstrictc{\mathcal G}(\pi_{\mathcal G},\SLpm{n})\bigr)
	\subset\Xref{n}(\mathcal G).
\end{aligned}
\]

For a nonempty semialgebraic set \(Y\), we write \(\dim Y\) for its semialgebraic dimension.
We adopt the convention \(\dim\varnothing=-\infty\).
The first main theorem shows that these spaces are semialgebraic and gives a common upper bound for their dimensions.
If \(\Xstr{n}(\mathcal G)\ne\varnothing\), both dimensions equal this bound.

\begin{theorem}\label{thm:main-dimension}
	Let \(\mathcal G\) be a \(2\)-perfect labeled truncation \(d\)-polytope, where \(d\ge3\).
	Assume that its Coxeter diagram is connected and that \(\mathcal G\) is not combinatorially equivalent to a \(d\)-simplex.
	Let \(f\) be the number of facets of \(\mathcal G\), and let \(e_2\) be the number of ridges labeled by \(2\).
	Let \(n\ge d+1\) be an integer, and set \(m=\min\{n,f\}\).
	Then the following statements hold:
	\begin{enumerate}[label=\textup{(\roman*)},itemsep=0.2ex,topsep=0.4ex]
		\item Both \(\Xref{n}(\mathcal G)\) and \(\Xstr{n}(\mathcal G)\) are semialgebraic.
		\item Their dimensions satisfy
		\[
		\dim\Xstr{n}(\mathcal G)
		\le
		\dim\Xref{n}(\mathcal G)
		\le
		(2m-d-2)f-m^2+\frac{d(d+1)}2+1-e_2.
		\]
		\item If \(\Xstr{n}(\mathcal G)\ne\varnothing\), then
		\[
		\dim\Xstr{n}(\mathcal G)
		=
		\dim\Xref{n}(\mathcal G)
		=
		(2m-d-2)f-m^2+\frac{d(d+1)}2+1-e_2.
		\]
	\end{enumerate}
	In particular, if \(\mathcal G\) is Coxeter realizable in \(\mathbb S^d\), then \(\Xstr{n}(\mathcal G)\ne\varnothing\), so \textup{(iii)} applies.
\end{theorem}

A representation is \emph{semisimple} if it is a direct sum of irreducible representations.
Under the hypotheses of Theorem~\ref{thm:main-dimension}, every character in \(\Xstr{n}(\mathcal G)\) admits a semisimple strict \(\mathcal G\)-compatible representative, unique up to \(\SLpm{n}\)-conjugacy (see the discussion preceding Definition~\ref{def:strict-character-rank-strata}).
For each \(x\in\Xstr{n}(\mathcal G)\), choose such a representative \(\rho\), and let \(C_\rho\) be an associated Cartan matrix.
The rank of \(C_\rho\) depends only on \(x\).
For each integer \(r\) with \(d+1\le r\le m\), set
\[
X^{\mathrm{ref},(r)}_{n,\mathrm{str}}(\mathcal G)
:=
\{x\in\Xstr{n}(\mathcal G)\mid\rank C_\rho=r\}.
\]

The second main theorem shows that these subspaces form a disjoint decomposition of \(\Xstr{n}(\mathcal G)\) and that each nonempty subspace is a smooth manifold.

\begin{theorem}\label{thm:strict-character-rank-strata}
	Let \(\mathcal G\) be a \(2\)-perfect labeled truncation \(d\)-polytope, where \(d\ge3\).
	Assume that its Coxeter diagram is connected and that \(\mathcal G\) is not combinatorially equivalent to a \(d\)-simplex.
	Let \(f\) be the number of facets of \(\mathcal G\), and let \(e_2\) be the number of ridges labeled by \(2\).
	Let \(n\ge d+1\) be an integer, and set \(m=\min\{n,f\}\).
	Then the following statements hold:
	\begin{enumerate}[label=\textup{(\roman*)},itemsep=0.2ex,topsep=0.4ex]
		\item The strict reflection character variety decomposes as
		\[
		\Xstr{n}(\mathcal G)
		=
		\bigsqcup_{r=d+1}^{m}X^{\mathrm{ref},(r)}_{n,\mathrm{str}}(\mathcal G).
		\]
		The set \(X^{\mathrm{ref},(m)}_{n,\mathrm{str}}(\mathcal G)\) is an open dense semialgebraic subset of \(\Xstr{n}(\mathcal G)\).
		Its relative boundary in \(\Xstr{n}(\mathcal G)\) is \(\bigsqcup_{r=d+1}^{m-1}X^{\mathrm{ref},(r)}_{n,\mathrm{str}}(\mathcal G)\).
		\item For each integer \(r\) with \(d+1\le r\le m\), the subspace \(X^{\mathrm{ref},(r)}_{n,\mathrm{str}}(\mathcal G)\), if nonempty, is a smooth manifold of dimension
		\[
		(2r-d-2)f-r^2+\frac{d(d+1)}2+1-e_2.
		\]
	\end{enumerate}
	In particular, if \(\mathcal G\) is Coxeter realizable in \(\mathbb S^d\), then every subspace in the decomposition is nonempty.
\end{theorem}

To prove the main theorems, we first show that, for \(d+1\le n\le f\), every character in \(\Xref{n}(\mathcal G)\) has a Vinberg-reduced representative.
Every nonempty fiber of \(\Xin\) restricted to the Vinberg-reduced subspace contains an \((n^2-1)\)-dimensional conjugacy orbit.
We therefore bound \(\dim\Xref{n}(\mathcal G)\) by estimating \(\dim\Homrefc{\mathcal G}(\pi_{\mathcal G},\SLpm{n})\) and subtracting \(n^2-1\).
A truncation sequence gives a facet ordering in which each facet and its earlier adjacent facets form a set of at most \(d+1\) pairwise adjacent facets.
We use this ordering and \(2\)-perfectness to prove that local scalar equations expressing the finite Coxeter relations have linearly independent differentials at \(\mathcal G\)-compatible representations.
The implicit function theorem then gives the required dimension bound and shows that the strict representation space, whenever nonempty, is a smooth manifold of dimension \((2n-d-2)f+\frac{d(d+1)}2-e_2\).

When \(d+1\le n\le f\) and \(\Homstrictc{\mathcal G}(\pi_{\mathcal G},\SLpm{n})\ne\varnothing\), a further dimension estimate shows that the representations of Cartan matrix rank \(n\) form an open dense subset of this space.
The conjugation action of \(\PGL{n}\) on this subset is free and proper, and its quotient is homeomorphic to \(X^{\mathrm{ref},(n)}_{n,\mathrm{str}}(\mathcal G)\).
The dimension of this quotient gives the matching lower bound, and its smooth structure transfers to the character subspace.
For \(n>f\), the Cartan matrix parametrization identifies the reflection character spaces with those for \(n=f\) by semialgebraic homeomorphisms.
Finally, for each integer \(r\) with \(d+1\le r\le m\), we apply the same quotient argument to strict \(\mathcal G\)-compatible representations in \(\SLpm{r}\) whose Cartan matrices have rank \(r\).
The block embedding \(\iota_r^n\) identifies their character space homeomorphically with \(X^{\mathrm{ref},(r)}_{n,\mathrm{str}}(\mathcal G)\) and transfers the resulting smooth structure.

The non-simplex assumption gives \(f\ge d+2\).
Thus the sets of at most \(d+1\) facets used in our independence arguments are proper subsets of \(S\), as required by Lemma~\ref{lem:adjacent-independence}.

The deformation space of Coxeter \(d\)-polytopes realizing \(\mathcal G\) is
\[
\Ccal(\mathcal G)
:=
\{\text{Coxeter realizations }(P,r_1,\ldots,r_f)
\text{ of }\mathcal G\text{ in }\mathbb S^d\}/\SLpm{d+1}.
\]
For this comparison, assume that \(\mathcal G\) and \(n\) satisfy the hypotheses of Theorem~\ref{thm:main-dimension} and that \(\Ccal(\mathcal G)\ne\varnothing\).
A truncation sequence gives \(e=df-\frac{d(d+1)}2\).
The results of Marquis~\cite[Theorem~3.16]{marquis2010espace} for \(d=3\) and of Choi-Lee-Marquis~\cite[Theorem~A]{choi2022deformation} for \(d\ge4\) give
\[
\dim\Ccal(\mathcal G)=e-e_2-d=df-\frac{d(d+3)}2-e_2.
\]

Since \(\Ccal(\mathcal G)\ne\varnothing\), Theorem~\ref{thm:main-dimension} gives
\[
\dim\Xref{n}(\mathcal G)
=
\dim\Xstr{n}(\mathcal G)
=
(2m-d-2)f-m^2+\frac{d(d+1)}2+1-e_2,
\]
where \(m=\min\{n,f\}\).
The difference between the dimension of either reflection character variety and that of \(\Ccal(\mathcal G)\) is therefore
\[
(m-d-1)(2f-m-d-1).
\]
Our results do not, in general, determine the number of connected components or the global topology of \(\Xref{n}(\mathcal G)\) and \(\Xstr{n}(\mathcal G)\).

An infinite discrete subgroup \(\Gamma\) of \(\PGL{n}\) is \emph{convex cocompact} in \(\mathbb{RP}^{n-1}\) if it acts properly discontinuously on a properly convex open set \(\Omega\subset\mathbb{RP}^{n-1}\) and cocompactly on a nonempty \(\Gamma\)-invariant convex subset \(C\subset\Omega\) that is closed in \(\Omega\) and whose closure in \(\mathbb{RP}^{n-1}\) contains all accumulation points of \(\Gamma\)-orbits in \(\Omega\)~\cite[Definition~1.1\textup{(ii)}]{danciger2025convex}.
When \(\pi_{\mathcal G}\) is infinite and word-hyperbolic, Danciger et al.~\cite[Corollary~1.11 and Proposition~2.7]{danciger2025convex} show that a \(\mathcal G\)-compatible representation \(\rho\) is strict if and only if its projective image is convex cocompact in \(\mathbb{RP}^{n-1}\).
Under the same hypotheses, strictness is also equivalent to \(\rho\) being \(P_1\)-Anosov~\cite[Corollary~1.18]{danciger2025convex}.

\begin{corollary}\label{cor:scc-reflection-characters}
	Let \(\mathcal G\) and \(n\) satisfy the hypotheses of Theorem~\ref{thm:main-dimension}, and suppose that \(\pi_{\mathcal G}\) is word-hyperbolic.
	Set \(m=\min\{n,f\}\).
	Then the following statements hold:
	\begin{enumerate}[label=\textup{(\roman*)},itemsep=0.2ex,topsep=0.4ex]
		\item The space \(\Xstr{n}(\mathcal G)\) is open in \(\Xn(\mathcal G)\) and consists precisely of the characters that admit a \(\mathcal G\)-compatible representative whose projective image is convex cocompact in \(\mathbb{RP}^{n-1}\).

		\item Let \(r\) be an integer with \(d+1\le r<m\), and let \(\bar\rho\in\Homstrictc{\mathcal G}(\pi_{\mathcal G},\SLpm{r})\).
		The projective image of the block extension \(\iota_r^n\circ\bar\rho\) is convex cocompact in \(\mathbb{RP}^{n-1}\).
		The character \(\Xin(\iota_r^n\circ\bar\rho)\) lies in the relative boundary of \(X^{\mathrm{ref},(m)}_{n,\mathrm{str}}(\mathcal G)\) in \(\Xstr{n}(\mathcal G)\).
		Every sufficiently small \(\mathcal G\)-compatible deformation of this block extension has convex cocompact projective image.
	\end{enumerate}
\end{corollary}

For odd \(n\ge3\), Danciger-Gu\'eritaud-Kassel~\cite[Proposition~1.7\textup{(1)}]{danciger2024projective} proved that Hitchin representations of closed orientable surface groups of genus at least \(2\) into \(\mathrm{PSL}(n,\mathbb R)\) have convex cocompact image in \(\mathbb{RP}^{n-1}\).

For a finite-dimensional real vector space \(V\), a subgroup \(\Gamma\) of \(\operatorname{PGL}(V)\) that is convex cocompact in \(\mathbb P(V)\) admits a lift \(\widehat{\Gamma}\) to \(\mathrm{SL}_{\pm}(V)\) preserving a properly convex cone in \(V\)~\cite[Remark~3.2]{danciger2024projective}.
The block embedding \(A\mapsto\diag(A,1)\) of \(\widehat{\Gamma}\) has convex cocompact projective image in \(\mathbb P(V\oplus\mathbb R)\)~\cite[Theorem~1.16\textup{(E)}]{danciger2024projective}.
Every sufficiently small deformation of the resulting representation \(\Gamma\to\operatorname{PGL}(V\oplus\mathbb R)\) has convex cocompact image~\cite[Theorem~1.16\textup{(D)}]{danciger2024projective}.

For an orderable Coxeter \(3\)-orbifold of normal type, Choi~\cite[Theorem~4]{choi2006deformation} proved that each nonempty restricted deformation space obtained by fixing the fundamental polyhedron up to projective equivalence is a smooth manifold (see Choi-Park~\cite[Theorem~4.7]{choi2025smoothness} for the formulation used here).
Choi-Park~\cite[Theorem~1.2]{choi2025smoothness} proved that the deformation space of an orderable labeled combinatorial \(3\)-polytope of normal type is a smooth manifold provided that its underlying polytope is not a cone over a polygon.
Both this space and the deformation space studied by Choi-Lee-Marquis~\cite[Definition~2.3]{choi2022deformation} parametrize Coxeter realizations of a fixed labeled combinatorial polytope.
We prescribe the Coxeter matrix of \(\mathcal G\) and allow the combinatorial type of the fundamental chamber in \(\mathbb S^{n-1}\) to vary.

For the Coxeter orbifold \(D^2(;n_1,n_2,n_3,n_4)\times\mathbb R\), where \(n_1,\ldots,n_4\) are integers at least \(3\), Bae~\cite[Remark~1.1 and Theorem~1.2]{bae2025convex} distinguished two cases according to whether the linear functionals defining the fundamental chamber are linearly independent.
He showed that the deformation space of convex projective structures with semisimple holonomy on this orbifold is homeomorphic to \(\mathbb R^3\times[4,\infty)^2\).
For a compact orientable \(2\)-orbifold \(O\) of negative Euler characteristic (possibly with boundary), Porti~\cite[Theorem~1.1]{porti2025deforming} studied block extensions of \(\mathbb C\)-irreducible representations \(\rho\colon\pi_1(O)\to\mathrm{SL}(k,\mathbb R)\) with \(k\ge2\).
He gave an explicit topological model for a neighborhood of \([\rho\oplus 1]\) in the \(\mathrm{SL}(k+1,\mathbb R)\)-character variety of \(\pi_1(O)\).

The remainder of the paper is organized as follows.
Section~\ref{sec:prelim} sets out our conventions for Coxeter systems, describes representations in terms of reflection pairs and Cartan matrices, and defines their characters using the real GIT quotient.
It also recalls the theorems of Marquis~\cite[Theorem~3.16]{marquis2010espace} and Choi-Lee-Marquis~\cite[Theorem~A]{choi2022deformation} on the deformation spaces of \(2\)-perfect labeled truncation polytopes.
Section~\ref{sec:Vinberg-reduced-representatives} constructs Vinberg-reduced representatives of reflection characters when \(n\le f\) and computes the centralizers of Vinberg-reduced representations when the Coxeter diagram is connected.
Section~\ref{sec:representation-space-dimensions} gives dimension bounds for the \(\mathcal G\)-compatible representation spaces and proves that the subspaces of strict representations are smooth manifolds whenever they are nonempty.
Section~\ref{sec:character-variety-dimensions} passes from representations to characters in the real GIT quotient and proves Theorem~\ref{thm:main-dimension}.
Section~\ref{sec:strict-rank-stratification} proves Theorem~\ref{thm:strict-character-rank-strata} on the decomposition of \(\Xstr{n}(\mathcal G)\) by Cartan matrix rank and Corollary~\ref{cor:scc-reflection-characters} on convex cocompactness.
Section~\ref{sec:structural-examples} gives three examples.
They illustrate the rank decomposition, exhibit a nonempty reflection character variety with an empty subspace of strict characters, and show that the \(2\)-perfectness assumption in Theorem~\ref{thm:main-dimension}\textup{(iii)} is essential.

\textbf{Acknowledgements.}
This work was carried out during the second author's doctoral studies at KAIST.
We thank Seungyeol Park, Juseop Lee, and Kwanho Kim for helpful discussions.
We are especially grateful to Seungyeol Park for the idea that led to Lemma~\ref{lem:block-embedded-strict} and for an earlier proof of the lemma.

\textbf{Use of AI tools.}
Yongho Seo used OpenAI's ChatGPT and Codex to assist with drafting, editing, and manuscript preparation.
The tools helped reformulate existing proofs of Lemmas~\ref{lem:one-facet-differential-criterion} and~\ref{lem:block-embedded-strict} in their present form.
The tools also provided initial ideas for the statement and proof of Proposition~\ref{prop:smooth-quotient-central-stabilizer}.
Yongho Seo has independently verified the mathematical arguments in these proofs.
We take full responsibility for the manuscript.

\section{Coxeter systems and reflection-generated representations}\label{sec:prelim}

In this section, we collect the definitions and preliminary results used throughout the paper.
Section~\ref{subsec:projective-coxeter-conventions} sets out our conventions for projective polytopes, following Choi-Lee-Marquis~\cite[Section~2.2]{choi2022deformation}, and defines the Coxeter system and Coxeter realizations of a labeled polytope.
Section~\ref{subsec:reflection-vinberg} states the Vinberg conditions and recalls that the spaces of \(\mathcal G\)-compatible representations are semialgebraic.
Section~\ref{subsec:reduced-strict-characters} defines Vinberg-reduced, dual-Vinberg-reduced, and strict representations and introduces the character varieties using the real GIT quotient.
Section~\ref{subsec:marquis-classification} recalls results on the topology and dimension of \(\Ccal(\mathcal G)\).

\subsection{Projective polytopes and associated Coxeter systems}\label{subsec:projective-coxeter-conventions}

The \emph{projective sphere} \(\mathbb S(V)\) of a finite-dimensional real vector space \(V\) is the space \((V\setminus\{0\})/\mathbb R_{>0}\) of positive rays.
If \(V=\mathbb R^{d+1}\), we write \(\mathbb S^d=\mathbb S(\mathbb R^{d+1})\).
For a subset \(A\subset V\), let \(\mathbb S(A)\) denote the image of \(A\setminus\{0\}\) in \(\mathbb S(V)\).
For \(n\ge1\), the group of projective automorphisms of \(\mathbb S^{n-1}\) is identified with
\[
\SLpm{n}=\{A\in \operatorname{GL}(n,\mathbb R)\mid \det A=\pm 1\}.
\]
Indeed, every class in \(\operatorname{GL}(n,\mathbb R)/\mathbb R_{>0}\) has a unique representative of determinant \(\pm1\).
The antipodal quotient of \(\mathbb S^{n-1}\) is \(\mathbb{RP}^{n-1}\).
The induced action on \(\mathbb{RP}^{n-1}\) gives the identification
\[
\SLpm{n}/\{\pm \Id_n\}\simeq\PGL{n}.
\]

The conjugation action of \(\SLpm{n}\) on \(\Hom(\pi_{\mathcal G},\SLpm{n})\) factors through \(\PGL{n}\), since the center acts trivially.
The induced \(\PGL{n}\)-action has the same orbits as the \(\SLpm{n}\)-action.

A \emph{projective reflection} of \(\mathbb S(V)\) is a projective automorphism of order two that fixes a projective hyperplane pointwise.
Equivalently, a projective reflection is induced by a linear map of the form
\[
R=\Id-\alpha\otimes v,
\qquad
\alpha\in V^*,\quad v\in V,\quad \alpha(v)=2.
\]
Here \((\alpha\otimes v)(x)=\alpha(x)v\) for \(x\in V\).
The hyperplane \(\mathbb S(\ker\alpha)\) is called the \emph{support} of \(R\).
The pair \((\alpha,v)\) is determined by \(R\) up to the rescaling
\[
(\alpha,v)\longmapsto(\lambda\alpha,\lambda^{-1}v),
\qquad \lambda\in\mathbb R^*.
\]
Fixing the ray \(\mathbb R_{>0}\alpha\) restricts the rescaling to \(\lambda>0\).

An \emph{affine chart} of \(\mathbb S(V)\) is a connected component of \(\mathbb S(V)\setminus\mathbb S(W)\), where \(W\) is a linear hyperplane of \(V\).
A subset \(\Omega\subset\mathbb S(V)\) is \emph{convex} if there exists a convex cone \(U\subset V\) such that \(\Omega=\mathbb S(U)\).
It is \emph{properly convex} if it is convex and its closure in \(\mathbb S(V)\) lies in an affine chart.
Equivalently, \(\Omega\) is convex and bounded in some affine chart.

For a nonzero linear form \(\alpha\in V^*\), the associated closed \emph{projective half-space} is
\[
H^-_\alpha
=
\mathbb S\bigl(\{x\in V\mid\alpha(x)\le0\}\bigr).
\]
A presentation as a finite intersection of half-spaces is \emph{irredundant} if deleting any half-space changes the intersection.
A \emph{projective polytope} is a properly convex subset \(P\subset\mathbb S(V)\) with nonempty interior in \(\mathbb S(V)\) that admits an irredundant presentation
\[
P=\bigcap_{i=1}^N H^-_{\alpha_i},
\qquad \alpha_1,\ldots,\alpha_N\in V^*\setminus\{0\}.
\]

In an affine chart containing \(P\), the faces of \(P\) are the usual faces of the corresponding compact affine polytope.
The face lattice \(\mathcal F(P)\) of \(P\) is the set of faces of \(P\) partially ordered by inclusion.
Two convex \(d\)-polytopes \(P\) and \(P'\) are combinatorially equivalent if their face lattices are isomorphic.

The empty face has rank \(-1\), and each nonempty face has rank equal to its dimension.
For \(d\ge2\), a \emph{combinatorial \(d\)-polytope} \(\mathcal G\) is a combinatorial equivalence class of convex \(d\)-polytopes.
A \emph{facet} of a \(d\)-polytope is a face of codimension \(1\), and a \emph{ridge} is a face of codimension \(2\).
An \emph{edge} is a \(1\)-dimensional face, and a \emph{vertex} is a \(0\)-dimensional face.

A combinatorial \(d\)-polytope is a \emph{truncation \(d\)-polytope} if it is represented by a convex \(d\)-polytope obtained from a \(d\)-simplex by a finite sequence of vertex truncations~\cite[Section~4.1]{choi2022deformation}.
A combinatorial \(d\)-polytope is \emph{simple} if each vertex lies in exactly \(d\) facets.
Every truncation \(d\)-polytope is simple.

\begin{definition}\label{def:labeledCox}
	Let \(\mathcal G\) be a combinatorial \(d\)-polytope with \(f\) facets, and let \(\mathcal F(\mathcal G)\) denote its face poset.
	A \emph{marking} is a fixed indexing of the facets of \(\mathcal G\) by \(\{1,\ldots,f\}\).
	We denote the indexed facets by \(\mathcal G_1,\ldots,\mathcal G_f\).
	For \(F,F'\in\mathcal F(\mathcal G)\), we write \(F\cap F'\) for their meet.
	Two facets are \emph{adjacent} if their intersection is a ridge.
	We label each ridge \(\mathcal G_i\cap\mathcal G_j\) by an integer \(m_{ij}\ge2\).
	A \emph{labeled combinatorial polytope} is a combinatorial polytope with a marking and a labeling of its ridges as above.
	An isomorphism between two such objects is an isomorphism of ranked face posets preserving the facet indices and the Coxeter labels.
	We also write \(\mathcal G\) for the resulting labeled combinatorial polytope.
\end{definition}

We write \(e\) for the number of ridges of \(\mathcal G\) and \(e_2\) for the number of ridges labeled by \(2\).
The associated Coxeter matrix \(M=(m_{ij})_{1\le i,j\le f}\) has diagonal entries \(m_{ii}=1\).
For \(i\ne j\), the entry \(m_{ij}\) is the label of \(\mathcal G_i\cap\mathcal G_j\) when the two facets are adjacent, and is \(\infty\) otherwise.
The Coxeter group associated with \(\mathcal G\) is
\[
\pi_{\mathcal G}
=
\left\langle s_1,\ldots,s_f \,\middle|\,
s_i^2=1\text{ for all }i,\quad
(s_i s_j)^{m_{ij}}=1
\text{ for }1\le i<j\le f\text{ with }m_{ij}<\infty
\right\rangle,
\]
where \(s_i\) is the generator corresponding to the facet \(\mathcal G_i\).

The Coxeter diagram \(\Sigma(\mathcal G)\) has nodes \(1,\ldots,f\), corresponding to the facets \(\mathcal G_1,\ldots,\mathcal G_f\).
Two distinct nodes \(i,j\) are joined by an edge exactly when \(m_{ij}\in\{3,4,\ldots,\infty\}\).
Each such edge has label \(m_{ij}\), with the label \(3\) omitted.
We follow Humphreys~\cite[Section~5.1]{humphreys1990reflection} for the conventions on Coxeter groups and diagrams.
We call \(\pi_{\mathcal G}\) \emph{irreducible} if \(\Sigma(\mathcal G)\) is connected.

Let \(S\) be the facet set of \(\mathcal G\).
For \(I\subseteq S\), let \(W_I\) denote the subgroup of \(\pi_{\mathcal G}\) generated by the \(s_i\) with \(\mathcal G_i\in I\).
We call \(W_I\) \emph{spherical} if it is finite.
For each edge \(E\) of \(\mathcal G\), set
\[
I_E:=\{s\in S\mid E\subseteq s\}.
\]
The labeled polytope \(\mathcal G\) is \emph{\(2\)-perfect} if \(W_{I_E}\) is spherical for every edge \(E\)~\cite[Section~1.2]{choi2022deformation}.
If \(\mathcal G\) is simple, then \(|I_E|=d-1\).
For \(d=3\), \(\mathcal G\) is automatically \(2\)-perfect, since each edge \(E\) is a ridge and its label determines the finite dihedral group \(W_{I_E}\).
A labeled polytope \(\mathcal G\) with a Coxeter realization \(P\) is \(2\)-perfect if and only if \(P\setminus\Omega_P\) is contained in the vertex set of \(P\)~\cite[Theorem~2.2\textup{(5)}]{choi2022deformation}.

The following definition is adapted from Choi-Park~\cite[Definition~2.2]{choi2025smoothness}.

\begin{definition}\label{def:coxeter-realization}
	Let \(d\ge2\), and let \(\mathcal G\) be a labeled combinatorial \(d\)-polytope with \(f\) facets \(\mathcal G_1,\ldots,\mathcal G_f\).
	A Coxeter \(d\)-polytope realizing \(\mathcal G\) is a projective \(d\)-polytope \(P\subset\mathbb S^d\) together with linear reflections \(r_1,\ldots,r_f\in\SLpm{d+1}\) subject to the following conditions.

	\begin{enumerate}[label=(\roman*),itemsep=0.2ex,topsep=0.4ex]
		\item There is a poset isomorphism \(\mathcal F(\mathcal G)\to\mathcal F(P)\) sending each \(\mathcal G_i\) to a facet \(P_i\) of \(P\).
		\item For each \(j\), the support of \(r_j\) is the hyperplane of \(\mathbb S^d\) containing \(P_j\).
		\item If \(\mathcal G_i\) and \(\mathcal G_j\) are adjacent, then \(r_i r_j\) has order \(m_{ij}\) in \(\SLpm{d+1}\).
		\item The subgroup \(\Gamma_P\subset\SLpm{d+1}\) generated by the linear reflections \(r_1,\ldots,r_f\) satisfies
		\[
		\operatorname{Int}(P)\cap\gamma\cdot\operatorname{Int}(P)=\varnothing
		\quad\text{for }\gamma\in\Gamma_P\setminus\{1\}.
		\]
	\end{enumerate}
	Such a tuple \((P,r_1,\ldots,r_f)\) is called a \emph{Coxeter realization} of \(\mathcal G\) in \(\mathbb S^d\).
	We say that \(\mathcal G\) is \emph{Coxeter realizable} in \(\mathbb S^d\) if it admits such a realization.
\end{definition}

\subsection{Reflection pairs, Cartan matrices, and Vinberg conditions}\label{subsec:reflection-vinberg}

We index the Coxeter matrix of \(\mathcal G\) by its facet set \(S\) and write \(M=(m_{st})_{s,t\in S}\).
We use the same symbol for a facet \(s\in S\) and the corresponding standard generator of \(\pi_{\mathcal G}\).

For \(n\ge1\), a representation \(\rho:\pi_{\mathcal G}\to\SLpm{n}\) is \emph{reflection-generated} if \(\rho(s)\) is a projective reflection for every \(s\in S\).
A \emph{reflection pair} for \(\rho\) is a family \((\alpha_s,v_s)_{s\in S}\), with \(\alpha_s\in(\mathbb R^n)^*\) and \(v_s\in\mathbb R^n\), such that
\[
\rho(s)=\Id-\alpha_s\otimes v_s,
\qquad
\alpha_s(v_s)=2
\quad\text{for every }s\in S.
\]
We define the Cartan matrix and Cartan products of the reflection pair by
\[
C=(c_{st})_{s,t\in S},
\qquad
c_{st}=\alpha_s(v_t),
\qquad
\mu_{st}=c_{st}c_{ts}.
\]

Let \(\lambda_s\in\mathbb R^*\) for every \(s\in S\), and set \(D=\diag(\lambda_s)\).
Replacing each pair \((\alpha_s,v_s)\) by \((\lambda_s\alpha_s,\lambda_s^{-1}v_s)\) sends \(C\) to \(DCD^{-1}\).
The products \(\mu_{st}\) are therefore independent of the choice of reflection pair for \(\rho\).
Fixing the rays \(\mathbb R_{>0}\alpha_s\) for all \(s\in S\) restricts the remaining rescalings to \(\lambda_s>0\).
We call two Cartan matrices \(C\) and \(C'\) \emph{diagonally equivalent} if \(C'=DCD^{-1}\) for some positive diagonal matrix \(D\).

The \emph{fundamental cone} associated with the reflection pair \((\alpha_s,v_s)_{s\in S}\) is
\[
\widetilde{\Delta}
=
\{x\in\mathbb R^n\mid \alpha_s(x)\le0\text{ for all }s\in S\}.
\]
The \emph{fundamental chamber} is its projectivization \(\mathbb S(\widetilde{\Delta})\subset\mathbb S^{n-1}\).

\begin{definition}\label{def:vinberg-conditions}
	Let \(\mathcal G\) be a labeled combinatorial polytope with facet set \(S\) and associated Coxeter group \(\pi_{\mathcal G}\).
	Let \(M=M(\mathcal G)=(m_{st})_{s,t\in S}\) be its associated Coxeter matrix.
	Let \(n\ge1\) be an integer, and let \(\rho:\pi_{\mathcal G}\to\SLpm{n}\) be reflection-generated.
	Choose a reflection pair \((\alpha_s,v_s)_{s\in S}\) for \(\rho\).
	Write \(C=(c_{st})_{s,t\in S}\) and \(\widetilde{\Delta}\subset\mathbb R^n\) for the Cartan matrix and fundamental cone of this reflection pair, respectively.
	This reflection pair satisfies the \emph{Vinberg conditions} for \(M\) if all of the following hold.
	\begin{enumerate}[label=(V\arabic*),labelindent=\parindent,leftmargin=*,itemsep=0.2ex,topsep=0.4ex]
		\item \(c_{st}=0\) for all \(s\ne t\) with \(m_{st}=2\).
		\item \(c_{st}<0\) for all \(s\ne t\) with \(m_{st}\ne2\).
		\item \(c_{st}c_{ts}=4\cos^2(\pi/m_{st})\) for all \(s\ne t\) with \(2<m_{st}<\infty\).
		\item \(c_{st}c_{ts}\ge4\) for all \(s\ne t\) with \(m_{st}=\infty\).
		\item \(\widetilde{\Delta}\) has nonempty interior in \(\mathbb R^n\).
	\end{enumerate}

	A real matrix \(C=(c_{st})_{s,t\in S}\) with \(c_{ss}=2\) for all \(s\in S\) is \emph{weakly compatible with \(M\)} if it satisfies \textup{(V1)}--\textup{(V3)}.
	The matrix \(C\) is \emph{compatible with \(M\)} if it is weakly compatible with \(M\) and satisfies \textup{(V4)}~\cite[Definitions~3.1 and~3.7]{danciger2025convex}.

	We call \(\rho\) \emph{\(\mathcal G\)-compatible} if it admits a reflection pair satisfying \textup{(V1)}--\textup{(V5)}.
	Equivalently, \(\rho\) is a representation of \(\pi_{\mathcal G}\) as a reflection group in the sense of Danciger et al.~\cite[Section~1.2]{danciger2025convex}.
	We denote the set of \(\mathcal G\)-compatible representations by \(\Homrefc{\mathcal G}(\pi_{\mathcal G},\SLpm{n})\).
\end{definition}

Thus every \(\mathcal G\)-compatible representation \(\rho\) is faithful and discrete.
In particular, the product \(\rho(s)\rho(t)\) has order \(m_{st}\) for all distinct \(s,t\in S\).
The Vinberg conditions do not require the associated fundamental chamber to be a projective polytope with face poset isomorphic to \(\mathcal F(\mathcal G)\).

An irreducible Coxeter group is \emph{affine} if it is infinite and has an abelian subgroup of finite index.

\begin{lemma}[Danciger et al.~{\cite[Section~3.3]{danciger2025convex}}]\label{lem:automatic-v5}
	Let \(\mathcal G\) be a labeled combinatorial polytope with facet set \(S\), Coxeter matrix \(M(\mathcal G)\), and Coxeter group \(\pi_{\mathcal G}\).
	Assume that \(\pi_{\mathcal G}\) is irreducible and non-affine.
	Let \(n\ge1\) be an integer, and let \(\rho:\pi_{\mathcal G}\to\SLpm{n}\) be reflection-generated.
	Choose a reflection pair \((\alpha_s,v_s)_{s\in S}\) for \(\rho\).
	If the Cartan matrix of this reflection pair is compatible with \(M(\mathcal G)\), then the pair satisfies \textup{(V5)}.
\end{lemma}

\begin{remark}\label{rem:truncation-group-hypotheses}
	Let \(d\ge3\), and let \(\mathcal G\) be a labeled truncation \(d\)-polytope with facet set \(S\).
	Assume that its Coxeter diagram is connected and that \(\mathcal G\) is not combinatorially equivalent to a \(d\)-simplex.
	Then \(\pi_{\mathcal G}\) is infinite, irreducible, and non-affine.

	The connectedness of the Coxeter diagram implies that \(\pi_{\mathcal G}\) is irreducible.
	Choose a sequence of vertex truncations that produces the underlying combinatorial type of \(\mathcal G\) from a \(d\)-simplex.
	The sequence is nonempty, and each truncation adds one facet, so \(|S|\ge d+2\).
	The first truncation creates a new facet that is disjoint from the facet opposite the truncated vertex of the original \(d\)-simplex.
	At each subsequent truncation, both facets are intersected with the retained half-space, so they remain disjoint.
	Thus \(\mathcal G\) has disjoint facets \(s,t\in S\) with \(m_{st}=\infty\).
	The subgroup \(W_{\{s,t\}}\) is infinite dihedral, so \(\pi_{\mathcal G}\) is infinite.
	The only irreducible affine Coxeter diagram containing an edge labeled \(\infty\) is \(\widetilde A_1\), which has two nodes~\cite[Appendix~A, Table~2]{danciger2025convex}.
	Since the Coxeter diagram of \(\mathcal G\) has more than two nodes, \(\pi_{\mathcal G}\) is non-affine.
\end{remark}

A subset of \(\mathbb R^N\) is \emph{semialgebraic} if it is a finite union of sets defined by finitely many polynomial equations and strict inequalities with real coefficients.
A map between semialgebraic sets is \emph{semialgebraic} if its graph is semialgebraic~\cite[Definition~2.2.5]{bochnak1998real}.
For nonempty semialgebraic sets, we use the dimension defined in~\cite[Definition~2.8.1]{bochnak1998real}.
We set \(\dim\varnothing=-\infty\).

\begin{theorem}[Tarski-Seidenberg {\cite[Theorem~2.2.1 and Proposition~2.2.7]{bochnak1998real}}]\label{thm:tarski-seidenberg}
	Let \(a,b\ge1\) be integers, and let \(p:\mathbb R^{a+b}\to\mathbb R^a\) be the projection onto the first \(a\) coordinates.
	If \(E\subset\mathbb R^{a+b}\) is semialgebraic, then \(p(E)\) is semialgebraic.
	In particular, if \(A\subset\mathbb R^a\) is semialgebraic and \(f:A\to\mathbb R^b\) has a semialgebraic graph, then \(f(A)\) is semialgebraic.
\end{theorem}

Every representation of a Coxeter group into \(\mathrm{GL}(n,\mathbb R)\) takes values in \(\SLpm{n}\), since the image of each standard generator squares to the identity and hence has determinant \(\pm1\).

\begin{proposition}[Danciger et al.~{\cite[Section~1.2]{danciger2025convex}}]\label{prop:homref-semialgebraic}
	Let \(\mathcal G\) be a labeled combinatorial polytope.
	For every integer \(n\ge1\), the set \(\Homrefc{\mathcal G}(\pi_{\mathcal G},\SLpm{n})\) is semialgebraic.
\end{proposition}

\subsection{Vinberg reduction, strictness, and GIT quotients}\label{subsec:reduced-strict-characters}

\begin{definition}\label{def:strict-reflective}
	Let \(\mathcal G\) be a labeled combinatorial polytope with facet set \(S\), Coxeter matrix \(M=(m_{st})_{s,t\in S}\), and Coxeter group \(\pi_{\mathcal G}\).
	Let \(n\ge1\) be an integer, and let \(\rho\in\Homrefc{\mathcal G}(\pi_{\mathcal G},\SLpm{n})\).
	Choose a reflection pair \((\alpha_s,v_s)_{s\in S}\) for \(\rho\).
	Set
	\[
	V_\alpha(\rho)=\bigcap_{s\in S}\ker\alpha_s,
	\qquad
	V_v(\rho)=\Span\{v_s\mid s\in S\}.
	\]
	The subspace \(V_\alpha(\rho)\) consists of the vectors fixed by \(\rho(\pi_{\mathcal G})\), and \(V_v(\rho)\) is the span of the \((-1)\)-eigenspaces of the generating reflections.
	Thus both subspaces are independent of the choice of reflection pair.
	We call \(\rho\) \emph{Vinberg-reduced} if \(V_\alpha(\rho)=\{0\}\).
	We call \(\rho\) \emph{dual-Vinberg-reduced} if \(V_v(\rho)=\mathbb R^n\).
	Danciger et~al.~\cite[Definition~3.5]{danciger2025convex} use the terms \emph{reduced} and \emph{dual-reduced}, respectively.

	The representation \(\rho\) is \emph{strict} if one, and hence every, Cartan matrix \(C=(c_{st})_{s,t\in S}\) of \(\rho\) satisfies
	\[
	c_{st}c_{ts}>4
	\quad\text{for all }s,t\in S\text{ with }m_{st}=\infty.
	\]
	We denote the set of strict \(\mathcal G\)-compatible representations by \(\Homstrictc{\mathcal G}(\pi_{\mathcal G},\SLpm{n})\).
\end{definition}

The identity \(c_{st}c_{ts}=\operatorname{tr}(\rho(s)\rho(t))-n+4\) shows that strictness is defined by finitely many strict polynomial inequalities in the matrix entries of the generating reflections.
Thus \(\Homstrictc{\mathcal G}(\pi_{\mathcal G},\SLpm{n})\) is semialgebraic by Proposition~\ref{prop:homref-semialgebraic}.
This set is relatively open in \(\Homrefc{\mathcal G}(\pi_{\mathcal G},\SLpm{n})\), but need not be dense.
In Example~\ref{ex:non-coxeter-realizable-empty-strict}, \(\Homrefc{\mathcal G}(\pi_{\mathcal G},\SLpm{5})\) is nonempty, whereas \(\Homstrictc{\mathcal G}(\pi_{\mathcal G},\SLpm{5})\) is empty.

For each integer \(n\ge1\), we write
\[
\begin{aligned}
	\Xn(\mathcal G)
	&:=\Hom(\pi_{\mathcal G},\SLpm{n})\sslash\SLpm{n},\\
	\Xin
	&:\Hom(\pi_{\mathcal G},\SLpm{n})\longrightarrow\Xn(\mathcal G),
\end{aligned}
\]
for the real GIT quotient under conjugation and the corresponding quotient map, respectively.
We define the reflection character variety and its subspace of strict characters by
\[
\begin{aligned}
	\Xref{n}(\mathcal G)
	&:=\Xin\bigl(\Homrefc{\mathcal G}(\pi_{\mathcal G},\SLpm{n})\bigr),\\
	\Xstr{n}(\mathcal G)
	&:=\Xin\bigl(\Homstrictc{\mathcal G}(\pi_{\mathcal G},\SLpm{n})\bigr)
	\subset\Xref{n}(\mathcal G).
\end{aligned}
\]

Fix an integer \(n\ge1\).
Evaluation on the standard generators realizes \(\Hom(\pi_{\mathcal G},\SLpm{n})\) as a closed algebraic subset of \(M_n(\mathbb R)^S\) invariant under simultaneous conjugation.
Following Richardson-Slodowy~\cite[Sections~7.2--7.3]{richardson1990minimum}, we identify \(\Xn(\mathcal G)\) with the set of closed \(\SLpm{n}\)-orbits.
The map \(\Xin\) sends each representation to the unique closed orbit in its orbit closure, taken in the usual topology.
Consequently, two representations have the same image under \(\Xin\) if and only if their orbit closures intersect.
Each fiber contains a unique closed orbit.
With the quotient topology induced by \(\Xin\), \(\Xn(\mathcal G)\) is Hausdorff~\cite[Theorem~7.7]{richardson1990minimum}.

Let \(\Gamma\) be a group and \(W\) a finite-dimensional real vector space.
A representation \(\rho:\Gamma\to\operatorname{GL}(W)\) is \emph{semisimple} if \(W\), viewed as a \(\Gamma\)-module via \(\rho\), is a direct sum of irreducible \(\Gamma\)-modules.
Equivalently, every \(\rho(\Gamma)\)-invariant subspace of \(W\) admits a \(\rho(\Gamma)\)-invariant complement.

\begin{lemma}[Danciger et al.~{\cite[Lemma~3.26]{danciger2025convex}}]\label{lem:semisimple-reflective}
	Let \(\mathcal G\) be a labeled combinatorial polytope whose Coxeter diagram is connected.
	Let \(n\ge1\) be an integer, and let \(\rho\in\Homrefc{\mathcal G}(\pi_{\mathcal G},\SLpm{n})\).
	Then \(\rho\) is semisimple if and only if
	\[
	\mathbb R^n=V_\alpha(\rho)\oplus V_v(\rho).
	\]
\end{lemma}

\begin{theorem}\label{thm:Richardson}
	Let \(\mathcal G\) be a labeled combinatorial polytope with associated Coxeter group \(\pi_{\mathcal G}\).
	Let \(n\ge1\) be an integer, and let \(\rho\in\Hom(\pi_{\mathcal G},\SLpm{n})\).
	Then the \(\SLpm{n}\)-orbit of \(\rho\) is closed in the usual topology on \(\Hom(\pi_{\mathcal G},\SLpm{n})\) if and only if \(\rho\) is semisimple.
\end{theorem}

\begin{proof}
	A conjugating matrix \(g\in\operatorname{GL}(n,\mathbb R)\) may be replaced by \(|\det g|^{-1/n}g\in\SLpm{n}\) without changing its conjugation action.
	Hence the conjugation actions of \(\operatorname{GL}(n,\mathbb R)\) and \(\SLpm{n}\) have the same orbits.
	By Richardson~\cite[1.2.1, Lemma~11.2, and Theorem~11.4]{richardson1988conjugacy}, the \(\operatorname{GL}(n,\mathbb R)\)-orbit of \((\rho(s))_{s\in S}\) is closed if and only if \(\rho\) is semisimple.
\end{proof}

By Theorem~\ref{thm:Richardson}, the semisimple representations in each fiber of \(\Xin\) form its unique closed orbit.
Thus every character in \(\Xn(\mathcal G)\) admits a semisimple representative, unique up to \(\SLpm{n}\)-conjugacy.

We use the results of Richardson-Slodowy~\cite[Section~6.1 and Proposition~6.7]{richardson1990minimum} to realize the real GIT quotient as the image of a polynomial map.

\begin{lemma}\label{lem:polynomial-real-GIT-realization}
	Let \(\mathcal G\) be a labeled combinatorial polytope with associated Coxeter group \(\pi_{\mathcal G}\).
	For every integer \(n\ge1\), there exist an integer \(N_n\ge1\) and a conjugation-invariant polynomial map
	\[
	p_n:\Hom(\pi_{\mathcal G},\SLpm{n})\longrightarrow\mathbb R^{N_n}.
	\]
	Two representations have the same image under \(p_n\) if and only if their \(\SLpm{n}\)-orbit closures in the usual topology intersect.

	Moreover, \(p_n\) induces a homeomorphism
	\[
	\bar p_n:\Xn(\mathcal G)\longrightarrow
	p_n\bigl(\Hom(\pi_{\mathcal G},\SLpm{n})\bigr),
	\]
	where the target has the subspace topology from \(\mathbb R^{N_n}\).
\end{lemma}

\begin{proof}
	Fix an integer \(n\ge1\).
	Using the standard generators, regard \(\Hom(\pi_{\mathcal G},\SLpm{n})\) as a closed algebraic subset of \(M_n(\mathbb R)^S\).
	The orbits under conjugation by \(\operatorname{GL}(n,\mathbb R)\) and by \(\SLpm{n}\) coincide.
	Let \(F_1,\ldots,F_{N_n}\) generate the \(\mathbb R\)-algebra of \(\operatorname{GL}(n,\mathbb R)\)-invariant polynomials on \(M_n(\mathbb R)^S\).
	Over \(\mathbb C\), the \(F_j\) generate the algebra of \(\operatorname{GL}(n,\mathbb C)\)-invariant polynomials on \(M_n(\mathbb C)^S\)~\cite[proof of Proposition~6.7]{richardson1990minimum}.
	Set
	\[
	p_n=(F_1,\ldots,F_{N_n})
	\big|_{\Hom(\pi_{\mathcal G},\SLpm{n})}.
	\]

	Suppose \(p_n(\rho)=p_n(\sigma)\), and let \(\rho_0\) and \(\sigma_0\) be semisimple representatives of \(\Xin(\rho)\) and \(\Xin(\sigma)\), respectively.
	By invariance and continuity, \(p_n(\rho_0)=p_n(\sigma_0)\).
	Their \(\operatorname{GL}(n,\mathbb C)\)-orbits are Zariski closed by Birkes~\cite[Corollary~5.3]{birkes1971orbits}.
	The complex invariant generators take the same values on \(\rho_0\) and \(\sigma_0\), so these orbits coincide~\cite[Section~6.1]{richardson1990minimum}.
	Let \(X+iY\) be an invertible complex intertwiner between \(\rho_0\) and \(\sigma_0\), with \(X,Y\in M_n(\mathbb R)\).
	The intertwining equations are real and linear, so \(X+tY\) is a real intertwiner for every \(t\in\mathbb R\).
	The polynomial \(\det(X+tY)\) is nonzero at \(t=i\), so \(X+tY\) is invertible for some real \(t\).
	Hence \(\rho_0\) and \(\sigma_0\) lie in the same \(\SLpm{n}\)-orbit, so the orbit closures of \(\rho\) and \(\sigma\) intersect.

	The converse follows from invariance and continuity of \(p_n\).
	Thus \(p_n\) induces a bijection \(\bar p_n\), which is continuous by the definition of the quotient topology.
	By Richardson-Slodowy~\cite[Proposition~6.7]{richardson1990minimum}, \(p_n\) sends closed invariant subsets of the representation space to closed subsets of \(\mathbb R^{N_n}\).
	If \(A\subseteq\Xn(\mathcal G)\) is closed, then \(\Xin^{-1}(A)\) is closed and invariant, so \(\bar p_n(A)=p_n(\Xin^{-1}(A))\) is closed.
	Hence \(\bar p_n\) is a homeomorphism.
\end{proof}

Using \(\bar p_n\), we henceforth identify \(\Xn(\mathcal G)\) with
\[
p_n\bigl(\Hom(\pi_{\mathcal G},\SLpm{n})\bigr)\subset\mathbb R^{N_n},
\]
so that \(\Xin=p_n\).

Let \(A\) be a semialgebraic subset of \(\Hom(\pi_{\mathcal G},\SLpm{n})\).
Since \(p_n\) is polynomial, Theorem~\ref{thm:tarski-seidenberg} implies that \(\Xin(A)\) is semialgebraic.

\subsection{Deformation spaces of labeled truncation polytopes}\label{subsec:marquis-classification}

The deformation space \(\Ccal(\mathcal G)\) consists of the projective isomorphism classes of Coxeter \(d\)-polytopes realizing \(\mathcal G\)~\cite[Definition~2.3]{choi2022deformation}.
To describe its topology and dimension, we use a consequence of Marquis's classification for \(d=3\) and the corresponding result of Choi-Lee-Marquis for \(d\ge4\).

\begin{theorem}[Marquis~{\cite[Theorem~3.16]{marquis2010espace}}]\label{thm:marquis-truncation-three}
	Let \(\mathcal G\) be a labeled truncation \(3\)-polytope whose Coxeter diagram is connected.
	Assume that \(\mathcal G\) is not combinatorially equivalent to a \(3\)-simplex and that \(\Ccal(\mathcal G)\) is nonempty.
	Then \(\Ccal(\mathcal G)\) is a finite disjoint union of open cells of dimension \(e-e_2-3\).
\end{theorem}

A Coxeter group is \emph{large} if it has a finite-index subgroup with a nonabelian free quotient.
Every infinite irreducible non-affine Coxeter group is large~\cite[Section~2.1]{choi2022deformation}, so \(\pi_{\mathcal G}\) is irreducible and large under the hypotheses of Remark~\ref{rem:truncation-group-hypotheses}.

\begin{theorem}[Choi-Lee-Marquis~{\cite[Theorem~A]{choi2022deformation}}]\label{thm:CLM-irreducible}
	Let \(d\ge4\), and let \(\mathcal G\) be a \(2\)-perfect labeled truncation \(d\)-polytope whose Coxeter group \(\pi_{\mathcal G}\) is irreducible and large.
	Assume that \(\Ccal(\mathcal G)\) is nonempty.
	Then \(d\le9\), and \(\Ccal(\mathcal G)\) is a finite disjoint union of open cells of dimension \(e-e_2-d\).
\end{theorem}

The bound \(d\le9\) in Theorem~\ref{thm:CLM-irreducible} concerns Coxeter realizations of \(\mathcal G\) in \(\mathbb S^d\).
Theorems~\ref{thm:main-dimension} and~\ref{thm:strict-character-rank-strata} apply for every \(d\ge3\) under their stated hypotheses, including when \(\Ccal(\mathcal G)=\varnothing\).

\section{Reduction to Vinberg-reduced representations}\label{sec:Vinberg-reduced-representatives}

In this section, we construct Vinberg-reduced representatives of reflection characters for the dimension estimates in Section~\ref{sec:character-variety-dimensions}.
Section~\ref{subsec:Vinberg-reduced-and-nonVinberg-reduced} defines the Vinberg-reduced and non-Vinberg-reduced subspaces of \(\Homrefc{\mathcal G}(\pi_{\mathcal G},\SLpm{n})\) and their images in the real GIT quotient.
For \(n\ge2\), Section~\ref{subsec:block-embedded} shows that every non-Vinberg-reduced representation has the same character as a block-embedded representation.
Section~\ref{subsec:centralizers-Vinberg-reduced-representations} computes the centralizers of Vinberg-reduced representations when the Coxeter diagram is connected.
Section~\ref{subsec:reduction-Vinberg-reduced-representatives} constructs a Vinberg-reduced representative for every reflection character when \(\mathcal G\) has at least \(n\) facets.
For a character in \(\Xstr{n}(\mathcal G)\), the representative can then be chosen strict.

\subsection{The Vinberg-reduced and non-Vinberg-reduced subspaces}\label{subsec:Vinberg-reduced-and-nonVinberg-reduced}

\begin{definition}\label{def:Vr-and-nVr-parts}
	Let \(\mathcal G\) be a labeled combinatorial polytope, and let \(n\ge1\) be an integer.
	We define the \emph{Vinberg-reduced subspace} of \(\Homrefc{\mathcal G}(\pi_{\mathcal G},\SLpm{n})\) by
	\[
	\HomVrc{\mathcal G}(\pi_{\mathcal G},\SLpm{n})
	:=
	\left\{
	\rho\in\Homrefc{\mathcal G}(\pi_{\mathcal G},\SLpm{n})
	\;\middle|\;
	\rho \text{ is Vinberg-reduced}
	\right\}.
	\]
	We define the \emph{non-Vinberg-reduced subspace} by
	\[
	\HomnVrc{\mathcal G}(\pi_{\mathcal G},\SLpm{n})
	:=
	\Homrefc{\mathcal G}(\pi_{\mathcal G},\SLpm{n})
	\setminus
	\HomVrc{\mathcal G}(\pi_{\mathcal G},\SLpm{n}).
	\]
\end{definition}

\begin{lemma}\label{lem:nVr-fixed-line}
	Let \(\mathcal G\) be a labeled combinatorial polytope, and let \(\pi_{\mathcal G}\) be its associated Coxeter group.
	Let \(n\ge1\) be an integer, and let \(\rho\in\Homrefc{\mathcal G}(\pi_{\mathcal G},\SLpm{n})\).
	The following conditions are equivalent:
	\begin{enumerate}[label=\textup{(\roman*)},itemsep=0.2ex,topsep=0.4ex]
		\item \(\rho\in\HomnVrc{\mathcal G}(\pi_{\mathcal G},\SLpm{n})\),
		\item the common fixed subspace of \(\rho(\pi_{\mathcal G})\) in \(\mathbb R^n\) is nonzero.
	\end{enumerate}
\end{lemma}

\begin{proof}
	By Definitions~\ref{def:strict-reflective} and~\ref{def:Vr-and-nVr-parts}, condition~\textup{(i)} holds if and only if \(V_\alpha(\rho)\ne\{0\}\).
	This is condition~\textup{(ii)}, since \(V_\alpha(\rho)\) is the common fixed subspace of \(\rho(\pi_{\mathcal G})\).
\end{proof}

Recall that \(\Xin\) is the real GIT quotient map.
We define the \emph{Vinberg-reduced character subspace} and the \emph{non-Vinberg-reduced character subspace} of \(\Xref{n}(\mathcal G)\) by
\[
\begin{aligned}
	\XnVr(\mathcal G)
	&:= \Xin\bigl(\HomVrc{\mathcal G}(\pi_{\mathcal G},\SLpm{n})\bigr),\\
	\XnnVr(\mathcal G)
	&:= \Xin\bigl(\HomnVrc{\mathcal G}(\pi_{\mathcal G},\SLpm{n})\bigr).
\end{aligned}
\]
These two subspaces need not be disjoint.

\subsection{Block-embedded representations}\label{subsec:block-embedded}

\begin{definition}\label{def:block-embedded-part}
	Let \(\mathcal G\) be a labeled combinatorial polytope, and let \(n\ge2\) be an integer.
	Let \(\iota\) be the standard block embedding
	\[
	\iota:\SLpm{n-1}\hookrightarrow\SLpm{n},
	\qquad
	\iota(A)=\diag(A,1).
	\]
	We call a representation \(\rho\in\Homrefc{\mathcal G}(\pi_{\mathcal G},\SLpm{n})\) \emph{block-embedded} if it is \(\SLpm{n}\)-conjugate to \(\iota\circ\rho_0\) for some \(\rho_0\in\Homrefc{\mathcal G}(\pi_{\mathcal G},\SLpm{n-1})\).
	We denote the space of block-embedded representations by \(\Homblockc{\mathcal G}(\pi_{\mathcal G},\SLpm{n})\).
	We define the \emph{block-embedded character subspace} of \(\Xref{n}(\mathcal G)\) by
	\[
	\Xnblock(\mathcal G)
	:=
	\Xin\bigl(\Homblockc{\mathcal G}(\pi_{\mathcal G},\SLpm{n})\bigr).
	\]
\end{definition}

\begin{lemma}\label{lem:nonVinberg-reduced-quotient-degeneration}
	Let \(\mathcal G\) be a labeled combinatorial polytope, and let \(n\ge1\) be an integer.
	Let \(\rho\) be a \(\mathcal G\)-compatible representation on \(V=\mathbb R^n\).
	Let \(W\) be a nonzero linear subspace of \(V_\alpha(\rho)\), and let \(\bar\rho\) be the induced representation on \(V/W\).
	Then the following statements hold:
	\begin{enumerate}[label=\textup{(\roman*)},itemsep=0.2ex,topsep=0.4ex]
		\item In any basis of \(V/W\), \(\bar\rho\) is \(\mathcal G\)-compatible.
		Every reflection pair for \(\rho\) satisfying \textup{(V1)}--\textup{(V5)} induces such a pair for \(\bar\rho\) with the same Cartan matrix.
		Moreover, \(V_\alpha(\bar\rho)=V_\alpha(\rho)/W\).
		\item Under any splitting \(V\cong(V/W)\oplus W\), the representation \(\rho_W=\bar\rho\oplus\Id_W\) is \(\mathcal G\)-compatible and lies in the closure of the \(\SLpm{n}\)-orbit of \(\rho\).
	\end{enumerate}
\end{lemma}

\begin{proof}
	    Set \(k=\dim W\) and \(r=\dim(V/W)=n-k\).
	We first prove \textup{(i)}.
	Let \((\alpha_s,v_s)_{s\in S}\) be any reflection pair for \(\rho\) satisfying \textup{(V1)}--\textup{(V5)}.
	Since \(\alpha_s(v_s)=2\), the common kernel \(V_\alpha(\rho)\) is a proper subspace of \(V\).
	As \(W\) is a nonzero subspace of \(V_\alpha(\rho)\), we have \(1\le r\le n-1\).
	The representation \(\rho\) fixes \(W\) pointwise, so \(\det\bar\rho(\gamma)=\det\rho(\gamma)\) for every \(\gamma\in\pi_{\mathcal G}\).
	Therefore \(\bar\rho\) takes values in \(\SLpm{r}\).
	Let \(q:V\to V/W\) be the quotient map.
	Since each \(\alpha_s\) vanishes on \(W\), it induces a linear form \(\bar\alpha_s\) on \(V/W\) with \(\alpha_s=\bar\alpha_s\circ q\).
	Set \(\bar v_s=q(v_s)\).
	For \(s,t\in S\), we have
	\[
	\bar\rho(s)=\Id_{V/W}-\bar\alpha_s\otimes\bar v_s,
	\qquad
	\bar\alpha_s(\bar v_t)=\alpha_s(v_t).
	\]
	Thus \((\bar\alpha_s,\bar v_s)_{s\in S}\) is a reflection pair for \(\bar\rho\) with the same Cartan matrix, so it satisfies \textup{(V1)}--\textup{(V4)}.

	Let \(\widetilde\Delta\subset V\) and \(\widetilde\Delta_W\subset V/W\) be the fundamental cones of the original and quotient reflection pairs, respectively.
	The identities \(\alpha_s=\bar\alpha_s\circ q\) give \(\widetilde\Delta=q^{-1}(\widetilde\Delta_W)\).
	Since the original pair satisfies \textup{(V5)} and \(q\) is open, \(q(\operatorname{Int}\widetilde\Delta)\) is a nonempty open subset of \(V/W\) contained in \(\widetilde\Delta_W\).
	Thus the quotient pair satisfies \textup{(V5)}, and \(\bar\rho\) is \(\mathcal G\)-compatible in any basis of \(V/W\).
	The factorization of the \(\alpha_s\) also gives \(q^{-1}\bigl(V_\alpha(\bar\rho)\bigr)=V_\alpha(\rho)\), hence \(V_\alpha(\bar\rho)=V_\alpha(\rho)/W\).

	To prove \textup{(ii)}, fix a linear splitting of \(q\) and use it to identify \(V\) with \((V/W)\oplus W\).
	Since \(\rho\) fixes \(W\) pointwise and induces \(\bar\rho\) on \(V/W\), we have
	\[
	\rho(\gamma)=
	\begin{pmatrix}
		\bar\rho(\gamma)&0\\
		b(\gamma)&\Id_W
	\end{pmatrix},
	\qquad \gamma\in\pi_{\mathcal G},
	\]
	where \(b(\gamma):V/W\to W\) is linear.
	A reflection pair for \(\rho_W\) is \(\bigl(\bar\alpha_s\oplus0,(\bar v_s,0)\bigr)_{s\in S}\).
	It has the same Cartan matrix as the quotient pair, and its fundamental cone \(\widetilde\Delta_W\times W\) has nonempty interior by \textup{(i)}.
	Thus \(\rho_W\) is \(\mathcal G\)-compatible.

	    For \(t>0\), set \(g_t=t^k\Id_{V/W}\oplus t^{-r}\Id_W\).
	Then \(\det g_t=t^{kr}t^{-rk}=1\), so \(g_t\in\SLpm{n}\).
	For \(\gamma\in\pi_{\mathcal G}\),
	\[
	g_t\rho(\gamma)g_t^{-1}
	=
	\begin{pmatrix}
		\bar\rho(\gamma)&0\\
		t^{-n}b(\gamma)&\Id_W
	\end{pmatrix}
	\longrightarrow \rho_W(\gamma)
	\qquad(t\to\infty).
	\]
	Convergence on the standard generators gives \(\rho_W\in\overline{\SLpm{n}\cdot\rho}\).
\end{proof}

\begin{corollary}\label{cor:git-nVr-equals-block-embedded}
	Let \(\mathcal G\) be a labeled combinatorial polytope, and let \(n\ge2\) be an integer.
	Then
	\[
	\XnnVr(\mathcal G)=\Xnblock(\mathcal G).
	\]
\end{corollary}

\begin{proof}
	By Definition~\ref{def:block-embedded-part}, every block-embedded representation fixes a line pointwise.
	Lemma~\ref{lem:nVr-fixed-line} therefore gives \(\Xnblock(\mathcal G)\subseteq\XnnVr(\mathcal G)\).

	For the reverse inclusion, let \(\rho\in\HomnVrc{\mathcal G}(\pi_{\mathcal G},\SLpm{n})\), and choose a line \(\ell\subseteq V_\alpha(\rho)\).
	Choose bases of \(\mathbb R^n/\ell\) and \(\ell\), and fix a splitting of the quotient map.
	Lemma~\ref{lem:nonVinberg-reduced-quotient-degeneration}\textup{(i)} gives a \(\mathcal G\)-compatible quotient representation \(\bar\rho\) into \(\SLpm{n-1}\).
	By part~\textup{(ii)}, the split representation \(\rho_\ell\) lies in \(\overline{\SLpm{n}\cdot\rho}\).
	In the resulting basis of \(\mathbb R^n\), the matrices of \(\rho_\ell\) are those of \(\iota\circ\bar\rho\).
	Rescaling the change-of-basis matrix to have determinant \(\pm1\) shows that \(\rho_\ell\) is block-embedded.
	Hence \(\Xin(\rho)=\Xin(\rho_\ell)\in\Xnblock(\mathcal G)\).
\end{proof}

\begin{remark}
	Let \(\mathcal G\) be a labeled combinatorial polytope, and let \(n\ge2\) be an integer.
	Suppose that \(\rho\in\HomnVrc{\mathcal G}(\pi_{\mathcal G},\SLpm{n})\) is not block-embedded.
	The proof of Corollary~\ref{cor:git-nVr-equals-block-embedded} gives a block-embedded representation \(\sigma\) in the closure of the \(\SLpm{n}\)-orbit of \(\rho\).
	By Lemma~\ref{lem:nVr-fixed-line}, \(\sigma\) is also non-Vinberg-reduced.
	Since being block-embedded is invariant under conjugation, we have \(\sigma\notin\SLpm{n}\cdot\rho\).
	Thus the orbit of \(\rho\) is not closed in the non-Vinberg-reduced subspace, and the orbit space of this subspace is non-Hausdorff.
	Since the orbit closures of \(\rho\) and \(\sigma\) meet, we have \(\Xin(\rho)=\Xin(\sigma)\).
\end{remark}

\subsection{Centralizers of Vinberg-reduced representations}\label{subsec:centralizers-Vinberg-reduced-representations}

\begin{proposition}\label{prop:centralizer-is-center}
	Let \(\mathcal G\) be a labeled combinatorial polytope with facet set \(S\) and connected Coxeter diagram \(\Sigma(\mathcal G)\).
	Let \(n\ge1\) be an integer, and let \(\rho:\pi_{\mathcal G}\to\SLpm{n}\) be a Vinberg-reduced \(\mathcal G\)-compatible representation.
	Then the centralizer \(Z_{\SLpm{n}}(\rho(\pi_{\mathcal G}))\) is the center \(\{\pm\Id_n\}\) of \(\SLpm{n}\).
\end{proposition}

\begin{proof}
	Let \(V=\mathbb R^n\), and let \(g\in Z_{\SLpm{n}}(\rho(\pi_{\mathcal G}))\).
	Choose a reflection pair \((\alpha_s,v_s)_{s\in S}\) for \(\rho\) satisfying \textup{(V1)}--\textup{(V5)}.
	Fix \(s\in S\).
	Since \(g\) commutes with \(\rho(s)\), it preserves the \((-1)\)-eigenspace \(\mathbb R v_s\), so \(g(v_s)=\lambda_s v_s\) for some \(\lambda_s\in\mathbb R^*\).
	The commutation relation \(g(\alpha_s\otimes v_s)=(\alpha_s\otimes v_s)g\) then gives \(\alpha_s\circ g=\lambda_s\alpha_s\).

	For \(s,t\in S\), the identities above give
	\[
	\lambda_s\alpha_s(v_t)
	=
	\alpha_s(g(v_t))
	=
	\lambda_t\alpha_s(v_t).
	\]
	If \(s\) and \(t\) are joined by an edge of \(\Sigma(\mathcal G)\), then \(m_{st}\ne2\), so \(\alpha_s(v_t)\ne0\) by condition~\textup{(V2)}.
	Hence \(\lambda_s=\lambda_t\).
	Since \(\Sigma(\mathcal G)\) is connected, all the scalars \(\lambda_s\) have a common value \(\lambda\).

	Since \(\rho\) is Vinberg-reduced, the covectors \((\alpha_s)_{s\in S}\) span \(V^*\).
	The identities \(\alpha_s\circ g=\lambda\alpha_s\) for all \(s\in S\) therefore give \(g=\lambda\Id_V\).
	Since \(|\det g|=|\lambda|^n=1\), we have \(\lambda=\pm1\).
	The matrices \(\pm\Id_n\) lie in the center of \(\SLpm{n}\).
\end{proof}

\subsection{Vinberg-reduced representatives in the GIT quotient}\label{subsec:reduction-Vinberg-reduced-representatives}

\begin{lemma}\label{lem:Vinberg-reduced-extension-from-lower-dimension}
	Let \(\mathcal G\) be a labeled combinatorial polytope with facet set \(S\).
	Let \(r\) and \(n\) be integers with \(1\le r<n\le |S|\), and let \(\bar\rho:\pi_{\mathcal G}\to\SLpm{r}\) be a Vinberg-reduced \(\mathcal G\)-compatible representation.
	Let \((\bar\alpha_s,\bar v_s)_{s\in S}\) be a reflection pair for \(\bar\rho\) satisfying \textup{(V1)}--\textup{(V5)}.
	Then there exists a Vinberg-reduced \(\mathcal G\)-compatible representation \(\rho:\pi_{\mathcal G}\to\SLpm{n}\) such that \(\bar\rho\oplus\Id_{n-r}\) lies in the closure of the \(\SLpm{n}\)-orbit of \(\rho\).
	Moreover, \(\rho\) admits a reflection pair that satisfies \textup{(V1)}--\textup{(V5)} and has the same Cartan matrix as \((\bar\alpha_s,\bar v_s)_{s\in S}\).
\end{lemma}

\begin{proof}
	Set \(U=\mathbb R^r\) and \(W=\mathbb R^{n-r}\), and identify \(V=U\oplus W\) with \(\mathbb R^n\).
	Define the evaluation map
	\[
	A_0:U\longrightarrow\mathbb R^S,\qquad
	A_0(u)=(\bar\alpha_s(u))_{s\in S}.
	\]
	Its kernel is the common kernel of the \(\bar\alpha_s\), which is zero because \(\bar\rho\) is Vinberg-reduced.
	A complement to \(A_0(U)\) in \(\mathbb R^S\) has dimension \(|S|-r\), which is at least \(\dim W=n-r\).
	Extend \(A_0\) to an injective linear map \(A:V\to\mathbb R^S\) by sending a basis of \(W\) to linearly independent vectors in such a complement.
	Write \(A=(\alpha_s)_{s\in S}\), and for each \(s\in S\) set
	\[
	v_s=(\bar v_s,0),\qquad R_s=\Id_V-\alpha_s\otimes v_s.
	\]

	Since \(A\) extends \(A_0\), we have \(\alpha_s|_U=\bar\alpha_s\) for each \(s\in S\).
	As \(v_t=(\bar v_t,0)\), it follows that \(\alpha_s(v_t)=\bar\alpha_s(\bar v_t)\) for all \(s,t\in S\).
	Thus \((\alpha_s(v_t))_{s,t\in S}\) is the original Cartan matrix, so it satisfies \textup{(V1)}--\textup{(V4)} and is weakly compatible with \(M(\mathcal G)\).
	By Danciger et al.~\cite[Section~3.2]{danciger2025convex}, the \(R_s\) are linear reflections satisfying the Coxeter relations, so \(s\mapsto R_s\) defines a representation \(\rho\) of \(\pi_{\mathcal G}\).
	The image of \(\rho\) lies in \(\SLpm{n}\) because it is generated by the reflections \(R_s\), each of which has determinant \(-1\).

	By \textup{(V5)} for the original pair, choose \(u_0\) in the interior of its fundamental cone.
	Since each \(\bar\alpha_s\) is nonzero, we have \(\alpha_s(u_0,0)=\bar\alpha_s(u_0)<0\) for every \(s\in S\).
	These finitely many strict inequalities hold on an open neighborhood of \((u_0,0)\) in \(V\).
	This proves \textup{(V5)} for \((\alpha_s,v_s)_{s\in S}\), and hence \(\rho\) is \(\mathcal G\)-compatible.
	Since \(A\) is injective, the common kernel of the \(\alpha_s\) is \(\ker A=\{0\}\), so \(\rho\) is Vinberg-reduced.

	For \(t>0\), set \(g_t=t^{n-r}\Id_U\oplus t^{-r}\Id_W\), with \(\det g_t=1\).
	For \(s\in S\), the definition of \(R_s\) gives \(\rho(s)(u,w)=\bigl(\bar\rho(s)u-\alpha_s(0,w)\bar v_s,w\bigr)\).
	Conjugation by \(g_t\) multiplies the upper-right block by \(t^n\), so
	\[
	g_t\rho(s)g_t^{-1}
	=
	\begin{pmatrix}
		\bar\rho(s) & -t^n\bigl(\alpha_s|_W\otimes\bar v_s\bigr)\\
		0 & \Id_W
	\end{pmatrix}
	\longrightarrow
	\bar\rho(s)\oplus\Id_W
	\qquad(t\to0).
	\]
	Convergence on the standard generators gives \(\bar\rho\oplus\Id_{n-r}\in\overline{\SLpm{n}\cdot\rho}\).
\end{proof}

\begin{theorem}\label{thm:reduction-to-Vinberg-reduced-representatives}
	Let \(\mathcal G\) be a labeled combinatorial polytope with facet set \(S\).
	Let \(n\) be an integer with \(1\le n\le |S|\).
	Then every character in \(\Xref{n}(\mathcal G)\) admits a Vinberg-reduced \(\mathcal G\)-compatible representative.
	For a character in \(\Xstr{n}(\mathcal G)\), the representative can be chosen strict.
\end{theorem}

\begin{proof}
	Let \([\sigma]\in\Xref{n}(\mathcal G)\), and choose a \(\mathcal G\)-compatible representative \(\tau\), taking \(\tau\) strict when \([\sigma]\in\Xstr{n}(\mathcal G)\).
	If \(V_\alpha(\tau)=\{0\}\), then \(\tau\) itself is the required representative.
	Otherwise, \(W=V_\alpha(\tau)\) is a nonzero subspace contained in a reflection hyperplane, so \(r=n-\dim W\) satisfies \(1\le r<n\).
	After conjugating \(\tau\) by an orthogonal matrix in \(\SLpm{n}\), we may assume \(W=\{0\}\oplus\mathbb R^{n-r}\).
	Choose a reflection pair for \(\tau\) satisfying \textup{(V1)}--\textup{(V5)}, and let \(C\) be its Cartan matrix.
	Lemma~\ref{lem:nonVinberg-reduced-quotient-degeneration}\textup{(i)} gives a \(\mathcal G\)-compatible quotient representation \(\bar\tau\) on \(\mathbb R^r\), which is Vinberg-reduced because \(W=V_\alpha(\tau)\).
	The induced reflection pair satisfies \textup{(V1)}--\textup{(V5)} and has Cartan matrix \(C\).
	With the standard splitting, part~\textup{(ii)} gives \(\bar\tau\oplus\Id_{n-r}\in\overline{\SLpm{n}\cdot\tau}\).

	Since \(1\le r<n\le |S|\), apply Lemma~\ref{lem:Vinberg-reduced-extension-from-lower-dimension} to \(\bar\tau\) and its induced reflection pair.
	The resulting Vinberg-reduced \(\mathcal G\)-compatible representation \(\rho\) admits a reflection pair with Cartan matrix \(C\), and its \(\SLpm{n}\)-orbit closure contains \(\bar\tau\oplus\Id_{n-r}\).
	Thus the two orbit closures meet, so
	\[
	\Xin(\rho)=\Xin(\tau)=[\sigma].
	\]
	If \(\tau\) is strict, then \(C\) satisfies the inequalities defining strictness, so \(\rho\) is strict as well.
\end{proof}

\section{Dimensions of \texorpdfstring{\(\mathcal G\)}{G}-compatible representation spaces}\label{sec:representation-space-dimensions}

Let \(d\ge3\), and let \(\mathcal G\) be a \(2\)-perfect labeled truncation \(d\)-polytope with \(f\) facets.
Assume that \(\Sigma(\mathcal G)\) is connected and that \(\mathcal G\) is not combinatorially equivalent to a \(d\)-simplex.
For \(n\ge d+1\), we bound the dimension of the \(\mathcal G\)-compatible representation space on \(\mathbb R^n\) and show that its subspace of strict representations is a smooth manifold whenever nonempty.
Section~\ref{subsec:facet-ordering-finite-edge-equations} orders the facets by a truncation sequence and writes the finite Coxeter relations as local equations.
Section~\ref{subsec:successive-facet-dimension-smoothness} proves that these equations have linearly independent differentials at each \(\mathcal G\)-compatible representation and applies the implicit function theorem.
For \(d+2\le n\le f\), Section~\ref{subsec:rank-deficient-strict-representations} bounds the dimension of the subspace of strict representations with Cartan matrices of rank less than \(n\).

\subsection{Facet ordering and local equations for finite Coxeter relations}\label{subsec:facet-ordering-finite-edge-equations}

Fix a sequence of vertex truncations that produces the underlying combinatorial type of \(\mathcal G\) from a \(d\)-simplex.
Reindex the facets so that \(\mathcal G_1,\ldots,\mathcal G_{d+1}\) are inherited from the original simplex.
Index the remaining facets as \(\mathcal G_{d+2},\ldots,\mathcal G_f\) in the order in which the truncations introduce them.
Use this indexing for the Coxeter matrix entries \(m_{ij}\) and the standard generators \(s_i\) of \(\pi_{\mathcal G}\).

We write \(e\) for the number of ridges of \(\mathcal G\) and \(e_2\) for the number of ridges labeled by \(2\).
For \(j=1,\ldots,f\), set
\[
I_j:=\{i<j\mid \mathcal G_i\cap\mathcal G_j\text{ is a ridge}\},
\qquad
Q_j:=\{i\in I_j\mid m_{ij}=2\},
\]
and write \(k_j:=|I_j|\) and \(q_j:=|Q_j|\).
The facets of the original simplex are pairwise adjacent.
At each truncation, the new facet is adjacent precisely to the \(d\) pairwise adjacent facets incident to the truncated vertex.
For \(d\ge3\), two existing facets are adjacent after a vertex truncation if and only if they were adjacent before it.
Hence \(k_j=j-1\) for \(1\le j\le d+1\) and \(k_j=d\) for \(d+2\le j\le f\).
Each ridge contributes exactly once to \(\sum_{j=1}^f k_j\), at the larger of its two facet indices.
A ridge labeled by \(2\) contributes to \(\sum_{j=1}^f q_j\) at the same index.
Hence
\[
\sum_{j=1}^f k_j=e
\qquad\text{and}\qquad
\sum_{j=1}^f q_j=e_2.
\]
Substituting the values of \(k_j\) gives \(e=df-\frac{d(d+1)}2\).

\begin{lemma}\label{lem:adjacent-independence}
	Let \(d\ge3\), and let \(\mathcal G\) be a \(2\)-perfect labeled truncation \(d\)-polytope with facet set \(S\).
	Assume that \(\Sigma(\mathcal G)\) is connected.
	Let \(n\ge d+1\) be an integer, and let \(\rho:\pi_{\mathcal G}\to\SLpm{n}\) be a \(\mathcal G\)-compatible representation.
	Choose a reflection pair \((\alpha_s,v_s)_{s\in S}\) for \(\rho\), and let \(C=(c_{st})_{s,t\in S}\) be the Cartan matrix of this pair.
	Let \(I\subsetneq S\) consist of at most \(d+1\) pairwise adjacent facets.
	Then the columns of \(C\) indexed by \(I\) are linearly independent, as are its rows indexed by \(I\).
	In particular, the families \((\alpha_s)_{s\in I}\) and \((v_s)_{s\in I}\) are linearly independent.
\end{lemma}

\begin{proof}
	Changing the reflection pair rescales each covector and vector by a nonzero factor and conjugates \(C\) by an invertible diagonal matrix.
	These changes preserve the linear independence of the rows, columns, covectors, and vectors indexed by \(I\).
	Since \(\rho\) is \(\mathcal G\)-compatible, we may therefore choose a reflection pair satisfying \textup{(V1)}--\textup{(V5)}.

	Suppose \(Ca=0\) for a nonzero vector \(a=(a_s)_{s\in S}\) with \(a_s=0\) for \(s\notin I\).
	If the nonzero coefficients of \(a\) have the same sign, replace \(a\) by \(-a\) if necessary so that they are positive.
	Since \(I\subsetneq S\), some coefficient of \(a\) is zero.
	Connectedness of \(\Sigma(\mathcal G)\) gives an edge joining \(s\) to \(t\) with \(a_s=0<a_t\).
	By \textup{(V1)} and \textup{(V2)}, the off-diagonal entries in row \(s\) are nonpositive and \(c_{st}<0\).
	Hence
	\[
	(Ca)_s=\sum_{u\ne s}c_{su}a_u\le c_{st}a_t<0.
	\]
	This contradicts \(Ca=0\), so \(a\) has both positive and negative coefficients.

	Replace \(a\) by \(-a\) if necessary so that at most half of its nonzero coefficients are positive.
	Set \(I_+=\{s\in I:a_s>0\}\). Then
	\[
	1\le |I_+|\le\frac{|I|}{2}\le\frac{d+1}{2}\le d-1.
	\]
	Write \(a_{I_+}=(a_s)_{s\in I_+}>0\) and \(C_{I_+}=(c_{st})_{s,t\in I_+}\).
	For \(s\in I_+\) and \(t\notin I_+\), both \(c_{st}\) and \(a_t\) are nonpositive.
	Thus \(Ca=0\) gives
	\[
	(C_{I_+}a_{I_+})_s=-\sum_{t\notin I_+}c_{st}a_t\le0
	\qquad(s\in I_+).
	\]

	If all facets in \(I_+\) are inherited from the original simplex, their intersection there has dimension \(d-|I_+|\).
	Otherwise, consider the truncation that introduces the last facet in \(I_+\).
	Each earlier facet in \(I_+\) is adjacent to the new facet in \(\mathcal G\), so the two facets are adjacent when the new facet is introduced.
	The new facet is a \((d-1)\)-simplex because the polytope before the cut is simple.
	The \(|I_+|-1\) earlier facets meet this simplex in distinct facets, so the common face again has dimension \(d-|I_+|\).
	In either case, the common face contains an edge because \(d-|I_+|\ge1\).
	At each later cut, at least one endpoint of the edge lies in the interior of the retained half-space.
	The retained segment is therefore again an edge.
	Let \(E\) be the resulting common edge in \(\mathcal G\).
	Since \(I_+\subseteq I_E\) and \(\mathcal G\) is \(2\)-perfect, \(W_{I_+}\) is finite.

	The nonnegative matrix \(B_{I_+}=2\Id-C_{I_+}\) is block diagonal, with one block for each irreducible factor of \(W_{I_+}\).
	Since \(W_{I_+}\) is finite, the least eigenvalue of each block of \(C_{I_+}\) is positive~\cite[Definition~3.12 and Fact~3.17\textup{(1)}]{danciger2025convex}.
	By Perron-Frobenius, the spectral radius of the corresponding block of \(B_{I_+}\) is \(2\) minus this eigenvalue and hence less than \(2\).
	Since \(C_{I_+}a_{I_+}\le0\), we have \((B_{I_+}/2)a_{I_+}\ge a_{I_+}>0\).
	Iterating with the nonnegative matrix \(B_{I_+}/2\) gives \((B_{I_+}/2)^ka_{I_+}\ge a_{I_+}\) for every \(k\ge1\), contradicting \((B_{I_+}/2)^k\to0\).

	Hence the columns of \(C\) indexed by \(I\) are linearly independent.
	The column argument uses only \(c_{ss}=2\) and \textup{(V1)}--\textup{(V3)} as conditions on the matrix.
	These conditions also hold for \(C^T\), so the rows of \(C\) indexed by \(I\) are linearly independent.
	Applying all the \(\alpha_t\) to a relation among the \(v_s\) with \(s\in I\) gives a relation among the corresponding columns.
	Evaluating a relation among the \(\alpha_s\) with \(s\in I\) on all the \(v_t\) gives a relation among the corresponding rows.
	Both relations are therefore trivial.
\end{proof}

For \(n\ge2\), let \(\Rcal_n\) be the smooth semialgebraic manifold of projective reflections of \(\mathbb S(\mathbb R^n)\).
The pairs \((\alpha,v)\) with \(\alpha(v)=2\) form a manifold of dimension \(2n-1\).
Quotienting by the one-dimensional rescaling action gives \(\dim\Rcal_n=2n-2\).

For a tuple \((R_1,\ldots,R_f)\in(\Rcal_n)^f\), set \(c_{ij}=\alpha_i(v_j)\), where \((\alpha_i,v_i)\) represents \(R_i\).
For \(n\ge d+1\), let \(\Ecal_n(\mathcal G)\) consist of the tuples for which the matrix \((c_{ij})\) satisfies \textup{(V1)} and \textup{(V3)}.
Rescaling preserves the zero entries and the products \(c_{ij}c_{ji}\), so the definition is independent of the choice of representing pairs.

Evaluation on the standard generators identifies \(\Homrefc{\mathcal G}(\pi_{\mathcal G},\SLpm{n})\) with a semialgebraic subset of \(\Ecal_n(\mathcal G)\).
This subset consists of the tuples that admit reflection pairs satisfying \textup{(V2)} and \textup{(V4)}.
Indeed, \textup{(V1)}--\textup{(V3)} ensure that \(s_i\mapsto R_i\) defines a representation~\cite[Section~3.2]{danciger2025convex}.
Lemma~\ref{lem:automatic-v5} then gives \textup{(V5)} from \textup{(V4)}, since \(\pi_{\mathcal G}\) is irreducible and non-affine by Remark~\ref{rem:truncation-group-hypotheses}.
For the strict representation space, replace \textup{(V4)} with \(c_{ij}c_{ji}>4\) for all \(i<j\) with \(m_{ij}=\infty\).
Extend a reflection pair satisfying the Vinberg conditions at a strict tuple to a continuous local choice.
Both \textup{(V2)} and the strict product inequalities persist, so the strict representation space is open in \(\Ecal_n(\mathcal G)\).

Fix a point \(R=(R_1,\ldots,R_f)\in\Ecal_n(\mathcal G)\).
For each \(i\), choose a coordinate in which a reflection vector for \(R_i\) is nonzero.
The reflections with this property form an open semialgebraic neighborhood \(\mathcal U_i\subset\Rcal_n\) of \(R_i\).
Normalizing the chosen coordinate to \(1\) on \(\mathcal U_i\) gives a smooth semialgebraic section of the projection from representing pairs to \(\Rcal_n\).
Write \((\alpha_i,v_i)\) for the resulting pair of local functions.

For \(j=1,\ldots,f\), define the smooth semialgebraic map
\[
\Phi_j:
\mathcal U_1\times\cdots\times\mathcal U_j
\longrightarrow
\mathbb R^{k_j+q_j}
\]
with components
\[
\begin{aligned}
	&\alpha_i(v_j),\ \alpha_j(v_i)
	&& (i\in Q_j),\\
	&\alpha_i(v_j)\alpha_j(v_i)-4\cos^2(\pi/m_{ij})
	&& (i\in I_j\setminus Q_j).
\end{aligned}
\]

Regard each \(\Phi_j\) as a map on the full product via projection onto the first \(j\) factors, and define
\[
\Phi:=(\Phi_1,\ldots,\Phi_f):
\mathcal U_1\times\cdots\times\mathcal U_f
\longrightarrow
\mathbb R^{e+e_2}.
\]
The zero set of \(\Phi\) is \(\Ecal_n(\mathcal G)\cap(\mathcal U_1\times\cdots\times\mathcal U_f)\).

\subsection{Differential surjectivity, dimension, and smoothness}\label{subsec:successive-facet-dimension-smoothness}

\begin{lemma}\label{lem:one-facet-differential-criterion}
	Let \(d\ge3\), and let \(\mathcal G\) be a \(2\)-perfect labeled truncation \(d\)-polytope with \(f\) facets.
	Assume that \(\Sigma(\mathcal G)\) is connected.
	Let \(n\ge d+1\) be an integer, and let \(R=(R_1,\ldots,R_f)\in\Homrefc{\mathcal G}(\pi_{\mathcal G},\SLpm{n})\).
	Fix \(j\in\{1,\ldots,f\}\), and let \(\Phi_j\) be a local map near \(R\) constructed as in Section~\ref{subsec:facet-ordering-finite-edge-equations}.
	For \(1\le i\le j\), write \((\alpha_i,v_i)\) for the representing pair at \(R_i\) used to define \(\Phi_j\).
	The differential in the \(j\)-th variable
	\[
	D_{R_j}\Phi_j:T_{R_j}\Rcal_n\longrightarrow\mathbb R^{k_j+q_j}
	\]
	is surjective if \(\alpha_j\notin\Span\{\alpha_i\mid i\in I_j\}\) or \(v_j\notin\Span\{v_i\mid i\in I_j\}\).
\end{lemma}

\begin{proof}
Suppose first that \(v_j\notin\Span\{v_i\mid i\in I_j\}\).
Keep \((\alpha_i,v_i)\) fixed for every \(i<j\).
Choose an arbitrary vector in \(\mathbb R^{k_j+q_j}\) with components \(\eta^v_{ij},\eta^\alpha_{ij}\) for \(i\in Q_j\) and \(\eta^\times_{ij}\) for \(i\in I_j\setminus Q_j\).
We seek a tangent vector \((\dot\alpha_j,\dot v_j)\) to the space of representing pairs at \((\alpha_j,v_j)\).
Differentiating \(\alpha_j(v_j)=2\) gives the tangent condition \(\dot\alpha_j(v_j)+\alpha_j(\dot v_j)=0\).
The prescribed component derivatives give the equations
\[
\begin{aligned}
	\alpha_i(\dot v_j)&=\eta^v_{ij}
	&& (i\in Q_j),\\
	\dot\alpha_j(v_i)&=\eta^\alpha_{ij}
	&& (i\in Q_j),\\
	\alpha_j(v_i)\alpha_i(\dot v_j)
	+\alpha_i(v_j)\dot\alpha_j(v_i)&=\eta^\times_{ij}
	&& (i\in I_j\setminus Q_j).
\end{aligned}
\]

By the facet ordering, \(I_j\) indexes at most \(d\) pairwise adjacent facets, all preceding \(\mathcal G_j\).
Lemma~\ref{lem:adjacent-independence} therefore gives linear independence of \((\alpha_i)_{i\in Q_j}\) and of \((v_i)_{i\in I_j}\).
We may then choose \(\dot v_j\in\mathbb R^n\) satisfying \(\alpha_i(\dot v_j)=\eta^v_{ij}\) for all \(i\in Q_j\).
The assumption on \(v_j\) implies that \(v_j\) and the \(v_i\) for \(i\in I_j\) are linearly independent.
For \(i\in I_j\setminus Q_j\), \(\alpha_i(v_j)\ne0\) because the Cartan product is positive by \textup{(V3)}.
Define a linear functional \(\dot\alpha_j\) on the span of these independent vectors by
\[
\begin{aligned}
	\dot\alpha_j(v_i)&=\eta^\alpha_{ij}
	&& (i\in Q_j),\\
	\dot\alpha_j(v_i)&=
	\frac{\eta^\times_{ij}-\alpha_j(v_i)\alpha_i(\dot v_j)}
	{\alpha_i(v_j)}
	&& (i\in I_j\setminus Q_j),\\
	\dot\alpha_j(v_j)&=-\alpha_j(\dot v_j),
\end{aligned}
\]
and extend it linearly to \(\mathbb R^n\).
The last equality gives the tangent condition.
The other equalities, together with the choice of \(\dot v_j\), give the prescribed component derivatives.

The projection from representing pairs to \(\Rcal_n\) has surjective differential at \((\alpha_j,v_j)\) because the chosen local section is a right inverse.
The kernel of this differential therefore has dimension \((2n-1)-(2n-2)=1\).
The rescaling orbit projects to \(R_j\), so its tangent line \(\mathbb R(\alpha_j,-v_j)\) is the kernel.
For \(i\in Q_j\), both component expressions vanish along the rescaling orbit by \textup{(V1)}.
For \(i\in I_j\setminus Q_j\), the Cartan product is constant along the orbit.
Thus the derivatives of the component expressions vanish on the kernel.
Project \((\dot\alpha_j,\dot v_j)\) to \(T_{R_j}\Rcal_n\) and apply the differential of the chosen local section to the resulting vector.
The lift and the original tangent vector have the same projection, so their difference lies in the kernel.
The lift therefore gives the prescribed component derivatives.
By the chain rule, \(D_{R_j}\Phi_j\) maps the projected vector to the chosen target.

For \(\alpha_j\notin\Span\{\alpha_i\mid i\in I_j\}\), solve the same equations in the reverse order.
Choose \(\dot\alpha_j\) first with \(\dot\alpha_j(v_i)=\eta^\alpha_{ij}\) for \(i\in Q_j\), using the independence of these \(v_i\).
For \(i\in I_j\setminus Q_j\), \(\alpha_j(v_i)\ne0\) by \textup{(V3)}, so the product equation determines \(\alpha_i(\dot v_j)\).
The covectors \(\alpha_j\) and \(\alpha_i\) for \(i\in I_j\) are linearly independent by Lemma~\ref{lem:adjacent-independence} and the assumption on \(\alpha_j\).
Thus we can choose \(\dot v_j\) to satisfy the tangent condition and all the equations for the component derivatives.
The projection argument gives surjectivity.
\end{proof}

\begin{corollary}\label{cor:full-differential-surjective}
	Let \(d\ge3\), and let \(\mathcal G\) be a \(2\)-perfect labeled truncation \(d\)-polytope with \(f\) facets and \(e\) ridges, of which \(e_2\) are labeled by \(2\).
	Assume that \(\Sigma(\mathcal G)\) is connected and that \(\mathcal G\) is not combinatorially equivalent to a \(d\)-simplex.
	Let \(n\ge d+1\) be an integer, and let \(R=(R_1,\ldots,R_f)\in\Homrefc{\mathcal G}(\pi_{\mathcal G},\SLpm{n})\).
	Let \(\Phi=(\Phi_1,\ldots,\Phi_f)\) be a local map near \(R\) as in Section~\ref{subsec:facet-ordering-finite-edge-equations}.
	Then
	\[
	D_R\Phi:
	\bigoplus_{i=1}^f T_{R_i}\Rcal_n
	\longrightarrow\mathbb R^{e+e_2}
	\]
	is surjective.
\end{corollary}

\begin{proof}
	For each \(j\), the set \(I_j\cup\{j\}\) indexes at most \(d+1\) pairwise adjacent facets.
	Since \(\mathcal G\) is not combinatorially equivalent to a \(d\)-simplex, we have \(f\ge d+2\).
	Hence \(I_j\cup\{j\}\) is a proper subset of \(\{1,\ldots,f\}\).
	By Lemma~\ref{lem:adjacent-independence}, the reflection vectors indexed by \(I_j\cup\{j\}\) are linearly independent, so \(v_j\notin\Span\{v_i\mid i\in I_j\}\).
	Thus, for every \(j\), Lemma~\ref{lem:one-facet-differential-criterion} shows that \(D_{R_j}\Phi_j\) is surjective.

	Each \(\Phi_j\) depends only on the first \(j\) factors, so \(D_R\Phi\) is block lower triangular.
	To realize an arbitrary target \((\eta_1,\ldots,\eta_f)\in\bigoplus_{j=1}^f\mathbb R^{k_j+q_j}\), choose \(\xi_j\in T_{R_j}\Rcal_n\) successively so that
	\[
	D_{R_j}\Phi_j(\xi_j)
	=\eta_j-\sum_{i<j}D_{R_i}\Phi_j(\xi_i).
	\]
	At step \(j\), the right-hand side is known, so surjectivity of \(D_{R_j}\Phi_j\) gives a suitable \(\xi_j\).
	The resulting tuple maps to the prescribed target under \(D_R\Phi\).
\end{proof}

\begin{proposition}\label{prop:full-representation-upper-bound}
	Let \(d\ge3\), and let \(\mathcal G\) be a \(2\)-perfect labeled truncation \(d\)-polytope with \(f\) facets and \(e_2\) ridges labeled by \(2\).
	Assume that \(\Sigma(\mathcal G)\) is connected and that \(\mathcal G\) is not combinatorially equivalent to a \(d\)-simplex.
	Let \(n\ge d+1\) be an integer.
	Then
	\[
	\dim\Homrefc{\mathcal G}(\pi_{\mathcal G},\SLpm{n})
	\le (2n-d-2)f+\frac{d(d+1)}{2}-e_2.
	\]
\end{proposition}

\begin{proof}
	If \(\Homrefc{\mathcal G}(\pi_{\mathcal G},\SLpm{n})=\varnothing\), the bound follows from \(\dim\varnothing=-\infty\).
	Otherwise, let \(R\) be any point of this representation space, and let \(\Phi\) be a local map near \(R\) constructed as in Section~\ref{subsec:facet-ordering-finite-edge-equations}.
	By Corollary~\ref{cor:full-differential-surjective}, \(D_R\Phi\) is surjective.
	The implicit function theorem applied to the local equation \(\Phi=0\) shows that \(\Ecal_n(\mathcal G)\) is a smooth manifold of dimension \((2n-2)f-e-e_2\) near \(R\).

	Since the representation space is a semialgebraic subset of \(\Ecal_n(\mathcal G)\), its local dimension at \(R\) is at most \((2n-2)f-e-e_2\).
	The dimension of a nonempty semialgebraic set is the maximum of its local dimensions~\cite[Proposition~2.8.12]{bochnak1998real}.
	Since \(R\) was arbitrary, the same bound holds for the dimension of the representation space.
	Substituting \(e=df-\frac{d(d+1)}{2}\) gives the claimed bound.
\end{proof}

\begin{proposition}\label{prop:strict-representation-spaces-smooth}
	Let \(d\ge3\), and let \(\mathcal G\) be a \(2\)-perfect labeled truncation \(d\)-polytope with \(f\) facets and \(e_2\) ridges labeled by \(2\).
	Assume that \(\Sigma(\mathcal G)\) is connected and that \(\mathcal G\) is not combinatorially equivalent to a \(d\)-simplex.
	Let \(n\ge d+1\) be an integer.
	Then \(\Homstrictc{\mathcal G}(\pi_{\mathcal G},\SLpm{n})\) is either empty or a smooth semialgebraic manifold of dimension \((2n-d-2)f+\frac{d(d+1)}{2}-e_2\).

	Moreover, when \(n=d+1\), every representation in this space has Cartan matrix rank \(d+1\) and is both Vinberg-reduced and dual-Vinberg-reduced.
\end{proposition}

\begin{proof}
	Assume that \(\Homstrictc{\mathcal G}(\pi_{\mathcal G},\SLpm{n})\neq\varnothing\), and let \(R\) be any point of this space.
	The proof of Proposition~\ref{prop:full-representation-upper-bound} shows that \(\Ecal_n(\mathcal G)\) is a smooth manifold of dimension \((2n-2)f-e-e_2\) near \(R\).
	The strict representation space is an open semialgebraic subset of \(\Ecal_n(\mathcal G)\) by the description in Section~\ref{subsec:facet-ordering-finite-edge-equations}.
	It is therefore a smooth semialgebraic manifold of the same dimension.
	Substituting \(e=df-\frac{d(d+1)}{2}\) gives the claimed dimension.

	Now assume that \(n=d+1\).
	The \(d+1\) facets inherited from the original \(d\)-simplex are pairwise adjacent and form a proper subset of the facet set.
	By Lemma~\ref{lem:adjacent-independence}, the corresponding columns of the Cartan matrix are linearly independent.
	These independent columns give rank at least \(d+1\), and the Cartan matrix has rank at most \(d+1\) because it factors through \(\mathbb R^{d+1}\).
	The same lemma shows that the corresponding reflection covectors and reflection vectors form bases of \((\mathbb R^{d+1})^*\) and \(\mathbb R^{d+1}\), respectively.
	Thus \(R\) is both Vinberg-reduced and dual-Vinberg-reduced.
\end{proof}

\subsection{Strict representations with Cartan matrices of rank less than \texorpdfstring{\(n\)}{n}}\label{subsec:rank-deficient-strict-representations}

We use the following dimension inequalities for semialgebraic maps here and in Section~\ref{sec:character-variety-dimensions}.
Let \(\varphi:A\to B\) be a semialgebraic map with \(A\ne\varnothing\).
Then
\begin{equation}\label{eq:semialgebraic-fiber-upper-bound}
	\dim A
	\le
	\dim\varphi(A)
	+\sup_{y\in\varphi(A)}\dim\varphi^{-1}(y).
\end{equation}
Moreover, if \(\ell\) is a nonnegative integer and \(\dim\varphi^{-1}(y)\ge\ell\) for every \(y\in\varphi(A)\), then
\begin{equation}\label{eq:semialgebraic-fiber-lower-bound}
	\dim\varphi(A)\le\dim A-\ell.
\end{equation}
These inequalities follow by applying semialgebraic triviality to the continuous projection from the graph of \(\varphi\) onto \(\varphi(A)\) and using the dimension properties of semialgebraic sets~\cite[Proposition~2.8.5 and Theorems~2.8.8 and~9.3.2]{bochnak1998real}.

Let \(n\ge d+1\) be an integer, and let \(\rho\) be a \(\mathcal G\)-compatible representation on \(V=\mathbb R^n\).
Choose a reflection pair \((\alpha_s,v_s)_{s\in S}\) for \(\rho\), and let \(C_\rho\) be the associated Cartan matrix.
Changing the reflection pair conjugates \(C_\rho\) by an invertible diagonal matrix and hence preserves its rank.
Define linear maps
\[
\begin{aligned}
	A_\rho &: V\longrightarrow\mathbb R^f,
	& A_\rho(x)&=(\alpha_s(x))_s,\\
	B_\rho &: \mathbb R^f\longrightarrow V,
	& B_\rho((t_s)_s)&=\sum_s t_s v_s.
\end{aligned}
\]
Then \(C_\rho=A_\rho B_\rho\).

Since \(\dim V=n\), we have
\begin{equation}\label{eq:rank-Vinberg-reduced-dual}
	\begin{aligned}
		\rank C_\rho=n
		&\quad\Longleftrightarrow\quad
		A_\rho\text{ is injective and }B_\rho\text{ is surjective}\\
		&\quad\Longleftrightarrow\quad
		\rho\text{ is Vinberg-reduced and dual-Vinberg-reduced}.
	\end{aligned}
\end{equation}
The identities \(\ker A_\rho=V_\alpha(\rho)\) and \(\operatorname{im}B_\rho=V_v(\rho)\) give the second equivalence.

Set
\[
\begin{aligned}
	\Lcal_C
	&:=\{\rho\in\Homstrictc{\mathcal G}(\pi_{\mathcal G},\SLpm{n})
	\mid\rank C_\rho\le n-1\},\\
	\Lcal_\alpha
	&:=\{\rho\in\Homstrictc{\mathcal G}(\pi_{\mathcal G},\SLpm{n})
	\mid\dim\Span\{\alpha_s\mid s\in S\}\le n-1\},\\
	\Lcal_v
	&:=\{\rho\in\Homstrictc{\mathcal G}(\pi_{\mathcal G},\SLpm{n})
	\mid\dim\Span\{v_s\mid s\in S\}\le n-1\}.
\end{aligned}
\]

For \(\Lcal_\alpha\), consider the semialgebraic set of pairs \((\rho,u)\) where \(\rho\) is strict, \(u\in V\setminus\{0\}\), and \(\rho(s)u=u\) for every \(s\in S\).
Since \(\rho(s)u=u\) is equivalent to \(\alpha_s(u)=0\) for each \(s\in S\), the projection \((\rho,u)\mapsto\rho\) has image \(\Lcal_\alpha\).
Theorem~\ref{thm:tarski-seidenberg} therefore shows that \(\Lcal_\alpha\) is semialgebraic.
For a nonzero covector \(\lambda\in V^*\) and each \(s\in S\), \(\lambda\rho(s)=\lambda\) is equivalent to \(\lambda(v_s)=0\).
The same projection argument shows that \(\Lcal_v\) is semialgebraic.
Finally, the equivalence~\eqref{eq:rank-Vinberg-reduced-dual} gives \(\Lcal_C=\Lcal_\alpha\cup\Lcal_v\), so \(\Lcal_C\) is semialgebraic as well.

We also define the space of strict representations whose Cartan matrices have rank \(n\) by
\[
\Homranknc{\mathcal G}(\pi_{\mathcal G},\SLpm{n})
:=
\{\rho\in\Homstrictc{\mathcal G}(\pi_{\mathcal G},\SLpm{n})\mid\rank C_\rho=n\}.
\]
This space is the complement of \(\Lcal_C\) in the strict representation space and is therefore semialgebraic.

\begin{lemma}\label{lem:fixed-line-quotient-dimension}
	Let \(\mathcal G\) be a labeled combinatorial polytope with \(f\) facets and associated Coxeter group \(\pi_{\mathcal G}\).
	Let \(n\ge2\) be an integer, and let \(\Lcal_\alpha\) be the set of strict \(\mathcal G\)-compatible representations on \(\mathbb R^n\) that fix a nonzero vector.
	Then
	\[
	\dim\Lcal_\alpha
	\le
	\dim\Homstrictc{\mathcal G}(\pi_{\mathcal G},\SLpm{n-1})+f+(n-1).
	\]
\end{lemma}

\begin{proof}
If \(\Lcal_\alpha=\varnothing\), the assertion is immediate.
Otherwise, set \(V=\mathbb R^n\) and define
\[
\widetilde{\Lcal}_\alpha
:=
\{(\rho,[u])\in\Lcal_\alpha\times\mathbb S(V)
\mid\mathbb Ru\subset V_\alpha(\rho)\}.
\]
For the unit representative \(u\) of \([u]\), the equations \(\rho(s)u=u\) for all \(s\in S\) cut out \(\widetilde{\Lcal}_\alpha\) in the product of the strict representation space with \(\mathbb S(V)\).
Thus \(\widetilde{\Lcal}_\alpha\) is semialgebraic.
Every \(\rho\in\Lcal_\alpha\) fixes a nonzero vector, so the semialgebraic projection \((\rho,[u])\mapsto\rho\) is surjective.
Hence \(\dim\Lcal_\alpha\le\dim\widetilde{\Lcal}_\alpha\).

For \(i=1,\ldots,n\), set \(U_i:=\{[u]\in\mathbb S(V)\mid u_i\ne0\}\) and \(W_i:=\{x\in V\mid x_i=0\}\).
The sets \(U_1,\ldots,U_n\) form a semialgebraic open cover of \(\mathbb S(V)\).
For \([u]\in U_i\), the condition \(u_i\ne0\) gives \(V=W_i\oplus\mathbb Ru\).
Identify \(W_i\) with \(\mathbb R^{n-1}\) by omitting the \(i\)-th coordinate, and set
\[
\widetilde{\Lcal}_{\alpha,U_i}
:=
\{(\rho,[u])\in\widetilde{\Lcal}_\alpha\mid [u]\in U_i\}.
\]

Fix \(i\).
For each \((\rho,[u])\in\widetilde{\Lcal}_{\alpha,U_i}\), let \(\bar\rho_{[u]}\) be the representation induced by \(\rho\) on \(V/\mathbb Ru\).
The decomposition \(V=W_i\oplus\mathbb Ru\) identifies this quotient with \(W_i\cong\mathbb R^{n-1}\).
Lemma~\ref{lem:nonVinberg-reduced-quotient-degeneration}\textup{(i)} shows that \(\bar\rho_{[u]}\) is \(\mathcal G\)-compatible.
A reflection pair for \(\rho\) induces one for the quotient with the same Cartan matrix, so \(\bar\rho_{[u]}\) is strict.
Define
\[
F_i:\widetilde{\Lcal}_{\alpha,U_i}
\longrightarrow
U_i\times\Homstrictc{\mathcal G}(\pi_{\mathcal G},\SLpm{n-1}),
\qquad
F_i(\rho,[u])=([u],\bar\rho_{[u]}).
\]
The projection onto \(W_i\) along \(\mathbb Ru\) is given by \(x\mapsto x-(x_i/u_i)u\).
Applying this projection to \(\rho(s)x\) for \(x\in W_i\) gives \(\bar\rho_{[u]}(s)x\).
The entries of \(\bar\rho_{[u]}(s)\) are rational functions of the entries of \(\rho(s)\) and the coordinates of the unit representative \(u\), with denominator \(u_i\).
Since \(u_i\ne0\) on \(U_i\), \(F_i\) is semialgebraic.

Fix a nonempty fiber of \(F_i\) over \(([u],\bar\rho)\).
For each \((\rho',[u])\) in this fiber, \(\rho'\) fixes \(\mathbb Ru\) pointwise and induces \(\bar\rho\) on \(V/\mathbb Ru\).
Using the coordinate basis of \(W_i\) followed by the unit representative \(u\), we have
\[
\rho'(s)=
\begin{pmatrix}
	\bar\rho(s)&0\\
	b_s&1
\end{pmatrix},
\qquad
b_s\in(\mathbb R^{n-1})^*.
\]
The relation \(\rho'(s)^2=\Id\) gives \(b_s\bar\rho(s)=-b_s\).
Thus \(b_s\) lies in the one-dimensional left \((-1)\)-eigenspace of the reflection \(\bar\rho(s)\).
Extracting the rows \((b_s)_{s\in S}\) in this fixed basis gives a semialgebraic map from the fiber to the product of these eigenspaces.
The rows determine the generator tuple, so this map is injective.
Thus the fiber has dimension at most \(f\), the dimension of the product.

For each \(i\) with \(\widetilde{\Lcal}_{\alpha,U_i}\ne\varnothing\), Equation~\eqref{eq:semialgebraic-fiber-upper-bound} and the bound on the fiber dimensions give
\[
\begin{aligned}
	\dim\widetilde{\Lcal}_{\alpha,U_i}
	&\le \dim F_i(\widetilde{\Lcal}_{\alpha,U_i})+f\\
	&\le \dim U_i+\dim\Homstrictc{\mathcal G}(\pi_{\mathcal G},\SLpm{n-1})+f.
\end{aligned}
\]
Each \(U_i\) is a nonempty open subset of \(\mathbb S(V)\), so \(\dim U_i=n-1\).
Since \(\widetilde{\Lcal}_\alpha\) is the finite union of the sets \(\widetilde{\Lcal}_{\alpha,U_i}\), the same bound holds for \(\dim\widetilde{\Lcal}_\alpha\).
The inequality \(\dim\Lcal_\alpha\le\dim\widetilde{\Lcal}_\alpha\) now gives the asserted bound.
\end{proof}

\begin{proposition}\label{prop:rank-deficient}
	Let \(d\ge3\), and let \(\mathcal G\) be a \(2\)-perfect labeled truncation \(d\)-polytope with \(f\) facets and \(e_2\) ridges labeled by \(2\).
	Assume that \(\Sigma(\mathcal G)\) is connected and that \(\mathcal G\) is not combinatorially equivalent to a \(d\)-simplex.
	Let \(n\) be an integer with \(d+2\le n\le f\).
	Let \(\Lcal_C\) be the set of strict \(\mathcal G\)-compatible representations on \(\mathbb R^n\) whose Cartan matrices have rank less than \(n\).
	Then
	\[
	\dim\Lcal_C
	\le (2n-d-3)f+n+\frac{d(d+1)}{2}-1-e_2.
	\]
\end{proposition}

\begin{proof}
	The equivalence~\eqref{eq:rank-Vinberg-reduced-dual} gives \(\Lcal_C=\Lcal_\alpha\cup\Lcal_v\).
	Applying Lemma~\ref{lem:fixed-line-quotient-dimension} and then Proposition~\ref{prop:strict-representation-spaces-smooth} in dimension \(n-1\), we obtain
	\[
	\begin{aligned}
		\dim\Lcal_\alpha
		&\le \dim\Homstrictc{\mathcal G}(\pi_{\mathcal G},\SLpm{n-1})+f+(n-1)\\
		&\le (2n-d-3)f+n+\frac{d(d+1)}{2}-1-e_2.
	\end{aligned}
	\]

	In the standard coordinates, we write vectors as columns and covectors as rows.
	For a representation \(\rho:\pi_{\mathcal G}\to\SLpm{n}\), define its dual by \(\mathcal D(\rho)(\gamma)=\rho(\gamma)^{-T}\) for \(\gamma\in\pi_{\mathcal G}\).
	For each standard generator \(s\), the relation \(\rho(s)^2=\Id\) gives \(\mathcal D(\rho)(s)=\rho(s)^T\).
	Thus \(\mathcal D\) is a semialgebraic involution acting by transposition on generator tuples.

	Let \(\rho\) be strict and \(\mathcal G\)-compatible, and choose a reflection pair \((\alpha_s,v_s)_{s\in S}\) satisfying \textup{(V1)}--\textup{(V5)}.
	Then \((v_s^T,\alpha_s^T)_{s\in S}\) is a reflection pair for \(\mathcal D(\rho)\) with Cartan matrix \(C_\rho^T\).
	Since \(M(\mathcal G)\) is symmetric and transposition preserves Cartan products, \(C_\rho^T\) satisfies \textup{(V1)}--\textup{(V4)} and the strict inequalities.
	By Remark~\ref{rem:truncation-group-hypotheses}, \(\pi_{\mathcal G}\) is irreducible and non-affine.
	Lemma~\ref{lem:automatic-v5} therefore gives \textup{(V5)} for this reflection pair.
	Thus \(\mathcal D\) restricts to an involution of the strict representation space.
	Since the reflection vectors and covectors are interchanged, \(\mathcal D(\Lcal_v)=\Lcal_\alpha\).

	The semialgebraic bijection \(\mathcal D:\Lcal_v\to\Lcal_\alpha\) gives \(\dim\Lcal_v=\dim\Lcal_\alpha\).
	Since \(\Lcal_C=\Lcal_\alpha\cup\Lcal_v\), we have
	\(\dim\Lcal_C=\max\{\dim\Lcal_\alpha,\dim\Lcal_v\}=\dim\Lcal_\alpha\),
	which gives the asserted bound.
\end{proof}

\begin{corollary}\label{cor:full-rank-strict-dense}
	Let \(d\ge3\), and let \(\mathcal G\) be a \(2\)-perfect labeled truncation \(d\)-polytope with \(f\) facets.
	Assume that \(\Sigma(\mathcal G)\) is connected and that \(\mathcal G\) is not combinatorially equivalent to a \(d\)-simplex.
	Let \(n\) be an integer with \(d+1\le n\le f\).
	Assume moreover that the strict representation space on \(\mathbb R^n\) is nonempty.
	Then \(\Homranknc{\mathcal G}(\pi_{\mathcal G},\SLpm{n})\) is a nonempty open dense subset of this space.
\end{corollary}

\begin{proof}
	For \(n=d+1\), Proposition~\ref{prop:strict-representation-spaces-smooth} shows that every strict representation has Cartan matrix rank \(n\), so the two spaces coincide.

	Assume henceforth that \(n\ge d+2\).
	Proposition~\ref{prop:strict-representation-spaces-smooth} shows that the strict representation space is a smooth semialgebraic manifold of dimension \((2n-d-2)f+\frac{d(d+1)}{2}-e_2\).
	The upper bound in Proposition~\ref{prop:rank-deficient} is smaller than this dimension by \(f-n+1>0\).
	Thus \(\Lcal_C\) has empty interior in this manifold.
	The complement of \(\Lcal_C\) is dense in the nonempty strict representation space and hence is itself nonempty.

	By Section~\ref{subsec:facet-ordering-finite-edge-equations}, we can choose Cartan matrices continuously near each representation.
	A nonzero \(n\times n\) minor remains nonzero nearby, so the rank-\(n\) condition is open.
\end{proof}

\section{Dimensions of the reflection character variety and its subspace of strict characters}\label{sec:character-variety-dimensions}

In this section, we pass from representations to characters to prove Theorem~\ref{thm:main-dimension}.
When \(d+1\le n\le f\) and the strict representation space is nonempty, Section~\ref{subsec:strict-character-dimension} computes the dimension of the subspace of strict characters of Cartan matrix rank \(n\).
This gives a lower bound for \(\dim\Xstr{n}(\mathcal G)\).
Section~\ref{subsec:full-character-variety} uses Vinberg-reduced representatives to obtain an upper bound for \(\dim\Xref{n}(\mathcal G)\) when \(d+1\le n\le f\).
For \(n>f\), the Cartan matrix parametrization reduces both dimension estimates to the case \(n=f\).
We also show that \(\Xref{n}(\mathcal G)\) and \(\Xstr{n}(\mathcal G)\) are semialgebraic.

\subsection{Strict reflection characters of Cartan matrix rank \texorpdfstring{\(n\)}{n}}\label{subsec:strict-character-dimension}

\begin{proposition}\label{prop:strict-character-lower-bound}
	Let \(d\ge3\), and let \(\mathcal G\) be a \(2\)-perfect labeled truncation \(d\)-polytope with \(f\) facets and \(e_2\) ridges labeled by \(2\).
	Assume that \(\Sigma(\mathcal G)\) is connected and that \(\mathcal G\) is not combinatorially equivalent to a \(d\)-simplex.
	Let \(n\) be an integer with \(d+1\le n\le f\).
	Assume moreover that the strict representation space on \(\mathbb R^n\) is nonempty.
	Then the following statements hold:
	\begin{enumerate}[label=\textup{(\roman*)},itemsep=0.2ex,topsep=0.4ex]
		\item \(\Homranknc{\mathcal G}(\pi_{\mathcal G},\SLpm{n})\) is a nonempty open subset of the strict representation space and is invariant under conjugation by \(\SLpm{n}\).
		It is a smooth semialgebraic manifold of dimension
		\[
		(2n-d-2)f+\frac{d(d+1)}{2}-e_2.
		\]
		\item Every representation in \(\Homranknc{\mathcal G}(\pi_{\mathcal G},\SLpm{n})\) is semisimple and has centralizer \(\{\pm\Id_n\}\) in \(\SLpm{n}\).
		\item The image \(\Xin(\Homranknc{\mathcal G}(\pi_{\mathcal G},\SLpm{n}))\) is semialgebraic and has dimension
		\[
		(2n-d-2)f-n^2+\frac{d(d+1)}{2}+1-e_2.
		\]
	\end{enumerate}
	In particular, the dimension in \textup{(iii)} is a lower bound for \(\dim\Xstr{n}(\mathcal G)\).
\end{proposition}

\begin{proof}
By Proposition~\ref{prop:strict-representation-spaces-smooth}, the nonempty strict representation space is a smooth semialgebraic manifold of the dimension in \textup{(i)}.
Corollary~\ref{cor:full-rank-strict-dense} shows that \(\Homranknc{\mathcal G}(\pi_{\mathcal G},\SLpm{n})\) is a nonempty open subset of this manifold.
This subset is semialgebraic by Section~\ref{subsec:rank-deficient-strict-representations}.
Conjugation preserves \(\mathcal G\)-compatibility and the Cartan matrix, so it preserves this subset.
This proves \textup{(i)}.

Let \(\rho\in\Homranknc{\mathcal G}(\pi_{\mathcal G},\SLpm{n})\).
By Equation~\eqref{eq:rank-Vinberg-reduced-dual}, \(\rho\) is both Vinberg-reduced and dual-Vinberg-reduced.
Lemma~\ref{lem:semisimple-reflective} therefore shows that \(\rho\) is semisimple.
By Proposition~\ref{prop:centralizer-is-center}, the centralizer of \(\rho\) in \(\SLpm{n}\) is \(\{\pm\Id_n\}\).
This proves \textup{(ii)}.

Since \(\rho\) is semisimple, its conjugacy orbit is closed by Theorem~\ref{thm:Richardson}.
The centralizer is finite, so the orbit dimension formula gives
\[
\dim\bigl(\SLpm{n}\cdot\rho\bigr)
=
\dim\SLpm{n}-\dim Z_{\SLpm{n}}\bigl(\rho(\pi_{\mathcal G})\bigr)
=
n^2-1.
\]

Consider the restriction of \(\Xin\) to the semialgebraic set \(\Homranknc{\mathcal G}(\pi_{\mathcal G},\SLpm{n})\).
Theorem~\ref{thm:tarski-seidenberg} shows that the image of this restriction is semialgebraic, since \(\Xin\) is polynomial in the realization of Lemma~\ref{lem:polynomial-real-GIT-realization}.
The domain is a union of closed conjugacy orbits.
Each fiber of \(\Xin\) contains a unique closed orbit, so every nonempty fiber of the restriction is a single orbit of dimension \(n^2-1\).

Equations~\eqref{eq:semialgebraic-fiber-upper-bound} and~\eqref{eq:semialgebraic-fiber-lower-bound} therefore give
\[
\begin{aligned}
	\dim\Xin(\Homranknc{\mathcal G}(\pi_{\mathcal G},\SLpm{n}))
	&=\dim\Homranknc{\mathcal G}(\pi_{\mathcal G},\SLpm{n})-(n^2-1)\\
	&=(2n-d-2)f-n^2+\frac{d(d+1)}{2}+1-e_2.
\end{aligned}
\]
This proves \textup{(iii)}.
The asserted lower bound follows because the image lies in \(\Xstr{n}(\mathcal G)\).
\end{proof}

\subsection{The full reflection character variety}\label{subsec:full-character-variety}

\begin{proposition}\label{prop:full-character-upper-bound}
	Let \(d\ge3\), and let \(\mathcal G\) be a \(2\)-perfect labeled truncation \(d\)-polytope with \(f\) facets and \(e_2\) ridges labeled by \(2\).
	Assume that \(\Sigma(\mathcal G)\) is connected and that \(\mathcal G\) is not combinatorially equivalent to a \(d\)-simplex.
	Let \(n\ge d+1\) be an integer, and set \(m=\min\{n,f\}\).
	Then \(\Xref{n}(\mathcal G)\) is semialgebraic and
	\[
	\dim\Xref{n}(\mathcal G)\le (2m-d-2)f-m^2+\frac{d(d+1)}{2}+1-e_2.
	\]
\end{proposition}

\begin{proof}
By Proposition~\ref{prop:homref-semialgebraic}, \(\Homrefc{\mathcal G}(\pi_{\mathcal G},\SLpm{n})\) is semialgebraic.
Theorem~\ref{thm:tarski-seidenberg} shows that the image \(\Xref{n}(\mathcal G)\) is semialgebraic, since \(\Xin\) is polynomial in the realization of Lemma~\ref{lem:polynomial-real-GIT-realization}.

First suppose that \(n\le f\).
The projection argument for \(\Lcal_\alpha\) in Section~\ref{subsec:rank-deficient-strict-representations} shows that the non-Vinberg-reduced \(\mathcal G\)-compatible representations form a semialgebraic set.
Its complement in \(\Homrefc{\mathcal G}(\pi_{\mathcal G},\SLpm{n})\) is the Vinberg-reduced subspace, which is therefore semialgebraic.
Conjugation preserves \(\mathcal G\)-compatibility and the dimension of the common fixed subspace, so the Vinberg-reduced subspace is invariant under conjugation.

By Theorem~\ref{thm:reduction-to-Vinberg-reduced-representatives}, every character in \(\Xref{n}(\mathcal G)\) admits a Vinberg-reduced \(\mathcal G\)-compatible representative.
Hence \(\Xref{n}(\mathcal G)=\XnVr(\mathcal G)\).

If \(\HomVrc{\mathcal G}(\pi_{\mathcal G},\SLpm{n})=\varnothing\), then \(\Xref{n}(\mathcal G)=\varnothing\).
The bound then follows from \(\dim\varnothing=-\infty\).
Otherwise, let \(\rho\in\HomVrc{\mathcal G}(\pi_{\mathcal G},\SLpm{n})\).
By Proposition~\ref{prop:centralizer-is-center}, the centralizer of \(\rho\) in \(\SLpm{n}\) is \(\{\pm\Id_n\}\).
The centralizer is finite, so the orbit dimension formula gives
\[
\dim\bigl(\SLpm{n}\cdot\rho\bigr)
=
\dim\SLpm{n}-\dim Z_{\SLpm{n}}\bigl(\rho(\pi_{\mathcal G})\bigr)
=
n^2-1.
\]
The Vinberg-reduced subspace is a union of conjugacy orbits.
Since \(\Xin\) is constant on each orbit, every nonempty fiber of its restriction to this subspace contains an orbit of dimension \(n^2-1\).

Applying Equation~\eqref{eq:semialgebraic-fiber-lower-bound} with \(\ell=n^2-1\) to this restriction and using Proposition~\ref{prop:full-representation-upper-bound}, we obtain
\[
\begin{aligned}
	\dim\Xref{n}(\mathcal G)
	&\le \dim\HomVrc{\mathcal G}(\pi_{\mathcal G},\SLpm{n})-(n^2-1)\\
	&\le \dim\Homrefc{\mathcal G}(\pi_{\mathcal G},\SLpm{n})-(n^2-1)\\
	&\le (2n-d-2)f-n^2+\frac{d(d+1)}{2}+1-e_2.
\end{aligned}
\]
Since \(m=n\), the required bound follows.

Now suppose that \(n>f\).
By Remark~\ref{rem:truncation-group-hypotheses}, \(\pi_{\mathcal G}\) is infinite, irreducible, and non-affine.
The parametrization of Danciger et al.~\cite[Facts~3.28--3.29]{danciger2025convex} uses weakly compatible Cartan matrices of rank at most the ambient dimension.
Lemma~\ref{lem:automatic-v5} supplies \textup{(V5)} for the semisimple representations associated with compatible matrices.
Every Cartan matrix has size \(f\times f\), so its rank is at most \(f\).
Thus the reflection character spaces in dimensions \(f\) and \(n\) are both identified homeomorphically with the same space of diagonal equivalence classes of compatible Cartan matrices.

The standard block embedding \(\iota_f^n\) preserves Cartan matrices and hence induces the identity on this common parameter space.
The resulting homeomorphism is
\[
\begin{aligned}
	\bar\iota_f^n:X^{\mathrm{ref}}_f(\mathcal G)
	&\longrightarrow\Xref{n}(\mathcal G),\\
	\Xi_{f,\mathcal G}(\rho)
	&\longmapsto\Xin(\iota_f^n\circ\rho).
\end{aligned}
\]
Since strictness is determined by the Cartan products, \(\bar\iota_f^n\) restricts to a homeomorphism from \(X^{\mathrm{ref}}_{f,\mathrm{str}}(\mathcal G)\) onto \(\Xstr{n}(\mathcal G)\).

In the polynomial realizations of Lemma~\ref{lem:polynomial-real-GIT-realization}, consider the polynomial map
\[
\rho\longmapsto\bigl(p_f(\rho),p_n(\iota_f^n\circ\rho)\bigr).
\]
The images of \(\Homrefc{\mathcal G}(\pi_{\mathcal G},\SLpm{f})\) and \(\Homstrictc{\mathcal G}(\pi_{\mathcal G},\SLpm{f})\) under this map are the graphs of \(\bar\iota_f^n\) and its restriction to \(X^{\mathrm{ref}}_{f,\mathrm{str}}(\mathcal G)\), respectively.
Both representation spaces are semialgebraic, so these graphs are semialgebraic by Theorem~\ref{thm:tarski-seidenberg}.
Interchanging the two factors gives semialgebraic graphs for the inverse maps, so both homeomorphisms preserve semialgebraic dimension.

Applying the bound already proved in dimension \(f\), we obtain
\[
\dim\Xref{n}(\mathcal G)
=
\dim X^{\mathrm{ref}}_f(\mathcal G)
\le
f^2-(d+2)f+\frac{d(d+1)}{2}+1-e_2.
\]
Since \(m=f\), the required bound follows.
\end{proof}

The following lemma adapts the argument of Choi-Lee~\cite[Proposition~2.6]{choiLee2015projective} to the present setting.

\begin{lemma}\label{lem:block-embedded-strict}
	Let \(d\ge2\), and let \(\mathcal G\) be a simple labeled combinatorial \(d\)-polytope.
	Assume that \(\Sigma(\mathcal G)\) is connected.
	Let \(\rho_0:\pi_{\mathcal G}\longrightarrow\SLpm{d+1}\) be a \(\mathcal G\)-compatible representation induced by a Coxeter realization of \(\mathcal G\).
	Let \(n\ge d+1\) be an integer.
	Then both \(\rho_0\) and \(\iota_{d+1}^n\circ\rho_0\) are strict.
\end{lemma}

\begin{proof}
Let \(P\) be the polytope of the given Coxeter realization.
Choose a reflection pair \((\alpha_s,v_s)_{s\in S}\) for \(\rho_0\) satisfying \textup{(V1)}--\textup{(V5)}, with fundamental chamber \(P\), and write \(c_{st}=\alpha_s(v_t)\).
Suppose that \(m_{st}=\infty\) and \(c_{st}c_{ts}=4\) for some \(s,t\in S\).
Set \(u=-c_{st}v_s+2v_t\).
Direct calculation gives
\[
\alpha_s(u)=\alpha_t(u)=0,
\qquad \alpha_j(u)=-c_{st}c_{js}+2c_{jt}\le0\quad(j\notin\{s,t\}).
\]

Since \(m_{st}=\infty\), the facets indexed by \(s\) and \(t\) are nonadjacent and hence disjoint by the simplicity of \(\mathcal G\).
If \(u\ne0\), the preceding equalities and inequalities imply that \([u]\) lies in \(P\) and in both facets, a contradiction.
Consequently \(u=0\).

For every \(j\in S\setminus\{s,t\}\), we have \(-c_{st}c_{js}+2c_{jt}=0\).
Both terms are nonpositive, and \(-c_{st}>0\), so \(c_{js}=c_{jt}=0\).
By \textup{(V2)}, \(m_{js}=m_{jt}=2\), so no edge of \(\Sigma(\mathcal G)\) joins \(\{s,t\}\) to \(S\setminus\{s,t\}\).
Since \(|S|\ge d+1\ge3\), \(\Sigma(\mathcal G)\) is disconnected, contrary to the hypothesis.
By \textup{(V4)}, we conclude that \(c_{st}c_{ts}>4\) whenever \(m_{st}=\infty\), so \(\rho_0\) is strict.

Finally, \(\bigl(\alpha_s\oplus0,(v_s,0)\bigr)_{s\in S}\) is a reflection pair for \(\iota_{d+1}^n\circ\rho_0\) with the same Cartan matrix.
Its fundamental cone is the product of the original fundamental cone with \(\mathbb R^{n-d-1}\), so \textup{(V5)} also holds.
Hence \(\iota_{d+1}^n\circ\rho_0\) is strict.
\end{proof}

\begin{proof}[Proof of Theorem~\ref{thm:main-dimension}]
	Proposition~\ref{prop:full-character-upper-bound} shows that \(\Xref{n}(\mathcal G)\) is semialgebraic and that
	\[
	\dim\Xref{n}(\mathcal G)
	\le
	(2m-d-2)f-m^2+\frac{d(d+1)}2+1-e_2.
	\]
	The discussion following Definition~\ref{def:strict-reflective} shows that \(\Homstrictc{\mathcal G}(\pi_{\mathcal G},\SLpm{n})\) is semialgebraic.
	Since \(\Xin\) is polynomial, Theorem~\ref{thm:tarski-seidenberg} implies that \(\Xstr{n}(\mathcal G)\) is semialgebraic, proving statement \textup{(i)}.
	Since \(\Xstr{n}(\mathcal G)\subset\Xref{n}(\mathcal G)\), the displayed upper bound also proves statement \textup{(ii)}.

	Assume that \(\Xstr{n}(\mathcal G)\ne\varnothing\).
	If \(n\le f\), then \(m=n\), so Proposition~\ref{prop:strict-character-lower-bound} gives
	\[
	\dim\Xstr{n}(\mathcal G)
	\ge
	(2m-d-2)f-m^2+\frac{d(d+1)}2+1-e_2.
	\]
	If \(n>f\), the map \(\bar\iota_f^n\) from the proof of Proposition~\ref{prop:full-character-upper-bound} restricts to a semialgebraic homeomorphism from \(X^{\mathrm{ref}}_{f,\mathrm{str}}(\mathcal G)\) onto \(\Xstr{n}(\mathcal G)\).
	Hence \(X^{\mathrm{ref}}_{f,\mathrm{str}}(\mathcal G)\ne\varnothing\), and Proposition~\ref{prop:strict-character-lower-bound}, applied in dimension \(f=m\), gives
	\[
	\dim\Xstr{n}(\mathcal G)
	=
	\dim X^{\mathrm{ref}}_{f,\mathrm{str}}(\mathcal G)
	\ge
	f^2-(d+2)f+\frac{d(d+1)}2+1-e_2.
	\]
	The lower bounds in the two cases and the inequalities in \textup{(ii)} prove statement \textup{(iii)}.

	Finally, assume that \(\mathcal G\) is Coxeter realizable in \(\mathbb S^d\), and let \(\rho_0:\pi_{\mathcal G}\longrightarrow\SLpm{d+1}\) be induced by a Coxeter realization.
	By Lemma~\ref{lem:block-embedded-strict}, \(\iota_{d+1}^n\circ\rho_0\) is a strict \(\mathcal G\)-compatible representation.
	Thus \(\Xstr{n}(\mathcal G)\ne\varnothing\), so statement \textup{(iii)} applies.
\end{proof}

\section{Rank decomposition of the strict reflection character variety}\label{sec:strict-rank-stratification}

We decompose \(\Xstr{n}(\mathcal G)\) by Cartan matrix rank and prove Theorem~\ref{thm:strict-character-rank-strata} and Corollary~\ref{cor:scc-reflection-characters}.
Under the hypotheses of Theorem~\ref{thm:strict-character-rank-strata}, every character in \(\Xstr{n}(\mathcal G)\) admits a semisimple strict \(\mathcal G\)-compatible representative by Danciger et al.~\cite[Facts~3.28--3.29]{danciger2025convex} and Lemma~\ref{lem:automatic-v5}.

\begin{definition}\label{def:strict-character-rank-strata}
	Let \(d\ge3\), and let \(\mathcal G\) be a \(2\)-perfect labeled truncation \(d\)-polytope with \(f\) facets.
	Assume that \(\Sigma(\mathcal G)\) is connected and that \(\mathcal G\) is not combinatorially equivalent to a \(d\)-simplex.
	Let \(n\ge d+1\) be an integer, and set \(m=\min\{n,f\}\).
	For each \(x\in\Xstr{n}(\mathcal G)\), choose a semisimple strict \(\mathcal G\)-compatible representative \(\rho\), and let \(C_\rho\) be any Cartan matrix of \(\rho\).
	For each integer \(r\) with \(d+1\le r\le m\), define the \emph{subspace of Cartan matrix rank \(r\)} by
	\[
	X^{\mathrm{ref},(r)}_{n,\mathrm{str}}(\mathcal G)
	:=
	\{x\in\Xstr{n}(\mathcal G)\mid\rank C_\rho=r\}.
	\]
	Any two choices of \(\rho\) and \(C_\rho\) yield Cartan matrices related by invertible diagonal conjugation, which preserves rank, so this subspace is well defined.
\end{definition}

\begin{lemma}\label{lem:strict-character-rank-bounds}
	Let \(d\ge3\), and let \(\mathcal G\) be a \(2\)-perfect labeled truncation \(d\)-polytope with \(f\) facets.
	Assume that \(\Sigma(\mathcal G)\) is connected and that \(\mathcal G\) is not combinatorially equivalent to a \(d\)-simplex.
	Let \(n\ge d+1\) be an integer, and set \(m=\min\{n,f\}\).
	Then every Cartan matrix associated with a character in \(\Xstr{n}(\mathcal G)\) has rank between \(d+1\) and \(m\).
\end{lemma}

\begin{proof}
	Fix \(x\in\Xstr{n}(\mathcal G)\), and choose a semisimple strict \(\mathcal G\)-compatible representative \(\rho\) with Cartan matrix \(C_\rho\), as in Definition~\ref{def:strict-character-rank-strata}.
	The \(d+1\) facets inherited from the original \(d\)-simplex are pairwise adjacent and form a proper subset of the facet set.
	Lemma~\ref{lem:adjacent-independence} shows that the corresponding \(d+1\) columns of \(C_\rho\) are linearly independent.
	Hence \(\rank C_\rho\ge d+1\).
	The factorization \(C_\rho=A_\rho B_\rho\) from Section~\ref{subsec:rank-deficient-strict-representations} gives \(\rank C_\rho\le n\), and \(C_\rho\) has size \(f\times f\), so \(\rank C_\rho\le m\).
\end{proof}

\begin{proposition}\label{prop:smooth-quotient-central-stabilizer}
	Let \(\mathcal G\) be a labeled combinatorial polytope, and let \(\pi_{\mathcal G}\) be its associated Coxeter group.
	Let \(n\ge1\) be an integer.
	Let \(\mathcal M\) be an \(\SLpm{n}\)-invariant smooth semialgebraic submanifold of \(\Hom(\pi_{\mathcal G},\SLpm{n})\).
	Assume that every \(\rho\in\mathcal M\) is semisimple and satisfies \(Z_{\SLpm{n}}(\rho(\pi_{\mathcal G}))=\{\pm\Id_n\}\).
	Then \(\Xin(\mathcal M)\) is semialgebraic and, if nonempty, is a smooth manifold of dimension \(\dim\mathcal M-(n^2-1)\).
\end{proposition}

\begin{proof}
	In the realization by invariant polynomials from Lemma~\ref{lem:polynomial-real-GIT-realization}, the map \(\Xin\) is polynomial.
	By Theorem~\ref{thm:tarski-seidenberg}, \(\Xin(\mathcal M)\) is semialgebraic.
	We may assume that \(\mathcal M\ne\varnothing\).

	Every representation in \(\mathcal M\) is irreducible over \(\mathbb R\).
	Indeed, if \(\rho\in\mathcal M\) had a nonzero proper invariant subspace, semisimplicity would give an invariant complement.
	The projection \(\Pi\) onto that subspace along its complement would commute with \(\rho\).
	The matrix \(\Id_n-2\Pi\) would act by \(-1\) on the subspace and by \(1\) on its complement.
	It would therefore be a nonscalar element of the centralizer in \(\SLpm{n}\), contrary to the hypothesis.

	The conjugation action of \(\PGL{n}\) on \(\mathcal M\) is proper.
	To verify this, let \(\rho_k\in\mathcal M\) and \(g_k\in\PGL{n}\) satisfy
	\[
	\rho_k\longrightarrow\rho,
	\qquad
	g_k\cdot\rho_k\longrightarrow\sigma,
	\qquad
	\rho,\sigma\in\mathcal M.
	\]
	Choose matrix representatives \(T_k\in\operatorname{GL}(n,\mathbb R)\) of \(g_k\), normalized by \(\lVert T_k\rVert=1\).
	After passing to a subsequence, we have \(T_k\to T\) with \(\lVert T\rVert=1\).
	For each generator \(s\), passing to the limit in \(T_k\rho_k(s)=(g_k\cdot\rho_k)(s)T_k\) gives
	\[
	T\rho(s)=\sigma(s)T.
	\]
	Thus \(\ker T\) is \(\rho\)-invariant.
	Since \(\rho\) is irreducible and \(T\ne0\), the matrix \(T\) is invertible.
	Hence \(g_k=[T_k]\to[T]\) in \(\PGL{n}\), proving properness by the sequential criterion.

	The \(\PGL{n}\)-action on \(\mathcal M\) is free.
	Indeed, every invertible matrix commuting with a representation in \(\mathcal M\) can be rescaled into \(\SLpm{n}\), where the centralizer hypothesis forces it to be scalar.
	The quotient manifold theorem~\cite[Theorem~9.18]{boumal2023introduction} gives \(\mathcal M/\PGL{n}\) a smooth manifold structure of dimension
	\[
	\dim\mathcal M-\dim\PGL{n}
	=
	\dim\mathcal M-(n^2-1).
	\]

	It remains to show that the natural map \(\mathcal M/\PGL{n}\to\Xin(\mathcal M)\) is a homeomorphism.
	Let \(\mathcal U\) be the subspace of representations in \(\Hom(\pi_{\mathcal G},\SLpm{n})\) that are irreducible over \(\mathbb R\).
	This subspace is open.
	For each integer \(k\) with \(1\le k<n\), the pairs \((\rho,L)\in\Hom(\pi_{\mathcal G},\SLpm{n})\times\operatorname{Gr}_k(\mathbb R^n)\) with \(L\) invariant under \(\rho\) form a closed incidence set.
	The Grassmannian is compact, so the image of this incidence set under projection to the representation space is closed.
	The set of reducible representations is the finite union of these images.

	Every irreducible representation is semisimple, so its conjugacy orbit is closed by Theorem~\ref{thm:Richardson}.
	We have \(\Xin^{-1}(\Xin(\mathcal U))=\mathcal U\).
	Otherwise, a reducible representation would have the same character as an irreducible one.
	The orbit closure of the reducible representation would then contain the closed orbit of the irreducible representation.
	This is impossible because the set of reducible representations is closed and invariant.
	Consequently, every nonempty fiber of \(\Xin|_{\mathcal U}\) is a single \(\PGL{n}\)-orbit.

	For any open subset \(O\subset\mathcal U\), the set \(\Xin^{-1}(\Xin(O))=\PGL{n}\cdot O\) is open in the full representation space.
	The quotient topology therefore makes \(\Xin(O)\) open in the real GIT quotient.
	Let \(O_{\mathcal M}\) be open in \(\mathcal M\).
	Choose an open subset \(O\subset\mathcal U\) with \(O_{\mathcal M}=\mathcal M\cap O\).
	Since the nonempty fibers of \(\Xin|_{\mathcal U}\) are conjugacy orbits and \(\mathcal M\) is invariant, we have
	\[
	\Xin(O_{\mathcal M})=\Xin(\mathcal M)\cap\Xin(O).
	\]
	Thus \(\Xin|_{\mathcal M}:\mathcal M\to\Xin(\mathcal M)\) is a continuous open surjection whose fibers are the conjugacy orbits.
	The induced bijection
	\[
	\mathcal M/\PGL{n}\longrightarrow\Xin(\mathcal M)
	\]
	is therefore a homeomorphism.
	Transporting the quotient smooth structure along this homeomorphism proves the assertion.
\end{proof}

\begin{proof}[Proof of Theorem~\ref{thm:strict-character-rank-strata}]
	Lemma~\ref{lem:strict-character-rank-bounds} gives the decomposition in \textup{(i)}.
	For each integer \(r\) with \(d+1\le r\le m\), set
	\[
	\mathcal M_r
	:=\{\rho\in\Homstrictc{\mathcal G}(\pi_{\mathcal G},\SLpm{r})\mid\rank C_\rho=r\}.
	\]
	This set is semialgebraic by Section~\ref{subsec:rank-deficient-strict-representations}.
	The block embedding \(\iota_r^n\) preserves Cartan matrices.
	By Remark~\ref{rem:truncation-group-hypotheses}, \(\pi_{\mathcal G}\) is infinite, irreducible, and non-affine.
	The parametrization of Danciger et al.~\cite[Facts~3.28--3.29]{danciger2025convex}, together with Lemma~\ref{lem:automatic-v5}, gives
	\begin{equation}\label{eq:rank-r-stratum-image}
		\Xin\bigl(\{\iota_r^n\circ\rho\mid\rho\in\mathcal M_r\}\bigr)
		=
		X^{\mathrm{ref},(r)}_{n,\mathrm{str}}(\mathcal G).
	\end{equation}
	It also gives a homeomorphism
	\[
	\begin{aligned}
		h_r^n:\Xi_{r,\mathcal G}(\mathcal M_r)
		&\longrightarrow X^{\mathrm{ref},(r)}_{n,\mathrm{str}}(\mathcal G),\\
		\Xi_{r,\mathcal G}(\rho)
		&\longmapsto\Xin(\iota_r^n\circ\rho).
	\end{aligned}
	\]
	Indeed, both spaces are parametrized by the same diagonal equivalence classes of compatible Cartan matrices of rank \(r\) satisfying the strict inequalities.
	In the polynomial realizations of Lemma~\ref{lem:polynomial-real-GIT-realization}, the graph of \(h_r^n\) is the image of \(\mathcal M_r\) under the polynomial map
	\[
	\rho\longmapsto\bigl(p_r(\rho),p_n(\iota_r^n\circ\rho)\bigr).
	\]
	By Theorem~\ref{thm:tarski-seidenberg}, this graph and its projection onto the second factor, \(X^{\mathrm{ref},(r)}_{n,\mathrm{str}}(\mathcal G)\), are semialgebraic.
	Interchanging the factors gives a semialgebraic graph for \((h_r^n)^{-1}\), so \(h_r^n\) and its inverse are semialgebraic.

	Under the parametrization by Cartan matrices, \(X^{\mathrm{ref},(m)}_{n,\mathrm{str}}(\mathcal G)\) is defined by the nonvanishing of at least one \(m\times m\) minor.
	It is therefore open in \(\Xstr{n}(\mathcal G)\).

	The density assertion is immediate if \(\Xstr{n}(\mathcal G)=\varnothing\).
	Otherwise, first suppose that \(n\le f\), so \(m=n\).
	Corollary~\ref{cor:full-rank-strict-dense} shows that \(\Homranknc{\mathcal G}(\pi_{\mathcal G},\SLpm{n})\) is dense in \(\Homstrictc{\mathcal G}(\pi_{\mathcal G},\SLpm{n})\).
	The continuity of \(\Xin\) and Equation~\eqref{eq:rank-r-stratum-image}, with \(r=n\), therefore imply that \(X^{\mathrm{ref},(n)}_{n,\mathrm{str}}(\mathcal G)\) is dense in \(\Xstr{n}(\mathcal G)\).

	Now suppose that \(n>f\), so \(m=f\).
	The homeomorphism from \(X^{\mathrm{ref}}_{f,\mathrm{str}}(\mathcal G)\) to \(\Xstr{n}(\mathcal G)\) constructed in the proof of Proposition~\ref{prop:full-character-upper-bound} preserves Cartan matrix rank.
	The density just proved in ambient dimension \(f\) therefore shows that \(X^{\mathrm{ref},(f)}_{n,\mathrm{str}}(\mathcal G)\) is dense in \(\Xstr{n}(\mathcal G)\).

	Since the top-rank subspace is open and dense, its relative boundary is its complement, namely the union of the subspaces of ranks \(d+1,\ldots,m-1\).
	For \(m=d+1\), this union is empty.
	This proves \textup{(i)}.

	Fix \(r\in\{d+1,\ldots,m\}\), and suppose that \(X^{\mathrm{ref},(r)}_{n,\mathrm{str}}(\mathcal G)\ne\varnothing\).
	Equation~\eqref{eq:rank-r-stratum-image} gives \(\mathcal M_r\ne\varnothing\).
	Since \(d+1\le r\le f\), Proposition~\ref{prop:strict-character-lower-bound} applies in dimension \(r\).
	Part \textup{(i)} of that proposition shows that \(\mathcal M_r\) is invariant under conjugation by \(\SLpm{r}\) and is a smooth semialgebraic manifold of dimension
	\[
	(2r-d-2)f+\frac{d(d+1)}2-e_2.
	\]
	Part \textup{(ii)} shows that every representation in this space is semisimple and has centralizer \(\{\pm\Id_r\}\).
	Proposition~\ref{prop:smooth-quotient-central-stabilizer}, applied in dimension \(r\), therefore shows that \(\Xi_{r,\mathcal G}(\mathcal M_r)\) is a smooth manifold of dimension
	\[
	(2r-d-2)f-r^2+\frac{d(d+1)}2+1-e_2.
	\]
	Transporting this smooth structure along \(h_r^n\) gives \(X^{\mathrm{ref},(r)}_{n,\mathrm{str}}(\mathcal G)\) a smooth manifold structure of the same dimension.
	This proves \textup{(ii)}.

	Finally, suppose that \(\mathcal G\) is Coxeter realizable in \(\mathbb S^d\).
	Let \(\rho_0:\pi_{\mathcal G}\to\SLpm{d+1}\) be induced by a Coxeter realization of \(\mathcal G\).
	Lemma~\ref{lem:block-embedded-strict} shows that \(\iota_{d+1}^r\circ\rho_0\) is strict for every integer \(r\) with \(d+1\le r\le m\).
	Thus the strict representation space on \(\mathbb R^r\) is nonempty for each such \(r\).
	Proposition~\ref{prop:strict-character-lower-bound}\textup{(i)}, applied in dimension \(r\), gives \(\mathcal M_r\ne\varnothing\).
	Equation~\eqref{eq:rank-r-stratum-image} then gives \(X^{\mathrm{ref},(r)}_{n,\mathrm{str}}(\mathcal G)\ne\varnothing\) for every such \(r\).
\end{proof}

\begin{remark}
	For \(d+1\le r<m\), when the subspaces of consecutive ranks \(r\) and \(r+1\) are both nonempty, their dimensions differ by \(2f-2r-1>0\).
	The homeomorphism induced by the block embedding transports the smooth quotient structure in dimension \(r\) to \(X^{\mathrm{ref},(r)}_{n,\mathrm{str}}(\mathcal G)\).
	Theorem~\ref{thm:strict-character-rank-strata} does not assert that \(\Xstr{n}(\mathcal G)\) is a manifold near a lower-rank subspace.
\end{remark}

\begin{proof}[Proof of Corollary~\ref{cor:scc-reflection-characters}]
	By Remark~\ref{rem:truncation-group-hypotheses}, \(\pi_{\mathcal G}\) is infinite, irreducible, and non-affine.
	Since \(\pi_{\mathcal G}\) is word-hyperbolic, Danciger et al.~\cite[Corollary~1.11 and Proposition~2.7]{danciger2025convex} show that a \(\mathcal G\)-compatible representation into \(\SLpm{n}\) is strict if and only if its projective image is convex cocompact in \(\mathbb{RP}^{n-1}\).
	This proves the characterization in \textup{(i)}.

	Every \(\operatorname{GL}(n,\mathbb R)\)-valued representation of \(\pi_{\mathcal G}\) takes values in \(\SLpm{n}\), since the standard generators are involutions.
	Moreover, \(\operatorname{GL}(n,\mathbb R)\)-conjugacy agrees with \(\SLpm{n}\)-conjugacy by the proof of Theorem~\ref{thm:Richardson}.
	The real GIT quotient of \(\Hom(\pi_{\mathcal G},\operatorname{GL}(n,\mathbb R))\) is therefore naturally homeomorphic to \(\Xn(\mathcal G)\).
	By Danciger et al.~\cite[Remark~3.4 and Fact~3.29]{danciger2025convex}, the characters represented by weakly compatible reflections form an open subset of \(\Xn(\mathcal G)\).
	Within this subset, the strict inequalities \(c_{st}c_{ts}>4\) for \(m_{st}=\infty\) define an open subset.
	Lemma~\ref{lem:automatic-v5} supplies \textup{(V5)}, so this open subset is precisely \(\Xstr{n}(\mathcal G)\).
	This completes the proof of \textup{(i)}.

	For \textup{(ii)}, let \(r\) be an integer with \(d+1\le r<m\), and let \(\bar\rho\in\Homstrictc{\mathcal G}(\pi_{\mathcal G},\SLpm{r})\).
	The block extension \(\iota_r^n\circ\bar\rho\) is \(\mathcal G\)-compatible and has the same Cartan matrix as \(\bar\rho\), so it is strict.
	By the equivalence proved above, its projective image is convex cocompact in \(\mathbb{RP}^{n-1}\).
	The Cartan matrix of \(\bar\rho\) factors through \(\mathbb R^r\), so its rank is at most \(r<m\).
	By Theorem~\ref{thm:strict-character-rank-strata}\textup{(i)}, the character \(\Xin(\iota_r^n\circ\bar\rho)\) therefore lies in the relative boundary of \(X^{\mathrm{ref},(m)}_{n,\mathrm{str}}(\mathcal G)\) in \(\Xstr{n}(\mathcal G)\).
	Finally, the strict inequalities persist under sufficiently small \(\mathcal G\)-compatible deformations.
	The same equivalence proves the stability assertion.
\end{proof}

\section{Examples}\label{sec:structural-examples}

We give three examples of reflection character varieties associated with labeled truncation polytopes.
In the first example, we use a normalized Cartan matrix to determine the connected components of the subspaces of Cartan matrix rank four and five in \(\Xstr{5}(\mathcal G)\).
The second example shows that \(\Xref{5}(\mathcal G)\) may be nonempty while \(\Xstr{5}(\mathcal G)\) is empty.
The third example shows that the \(2\)-perfectness assumption in Theorem~\ref{thm:main-dimension}\textup{(iii)} is essential.

\begin{definition}
	Let \(\mathcal G\) be a labeled combinatorial polytope with connected Coxeter diagram \(\Sigma(\mathcal G)\).
	Fix a spanning tree \(\mathcal T\) of \(\Sigma(\mathcal G)\), an orientation \(\overrightarrow{\mathcal T}\) of \(\mathcal T\), and a negative number \(b_{st}\) for each \((s,t)\in\overrightarrow{\mathcal T}\).
	By Greene~\cite[Lemma~2.2.1]{greene2013deformation}, each Cartan matrix \(C=(c_{st})\) compatible with \(M(\mathcal G)\) is diagonally equivalent to a matrix \(C_{\mathcal T}\) satisfying
	\[
	(C_{\mathcal T})_{st}=b_{st}
	\quad\text{for every }(s,t)\in\overrightarrow{\mathcal T}.
	\]
	This representative is unique, since connectedness of \(\mathcal T\) forces any positive diagonal matrix preserving the prescribed entries to be scalar.
	For every \((s,t)\in\overrightarrow{\mathcal T}\), the entry in the reverse direction is
	\[
	(C_{\mathcal T})_{ts}=\frac{\mu_{st}}{b_{st}},
	\]
	where \(\mu_{st}=c_{st}c_{ts}\) is invariant under diagonal equivalence.
	We call \(C_{\mathcal T}\) the \emph{\(\mathcal T\)-normalized Cartan matrix} with respect to these choices.
\end{definition}

\begin{example}\label{ex:prism-one-vertex-truncation}

	\begin{figure}[htbp]
		\centering
		\includegraphics[width=0.36\textwidth]{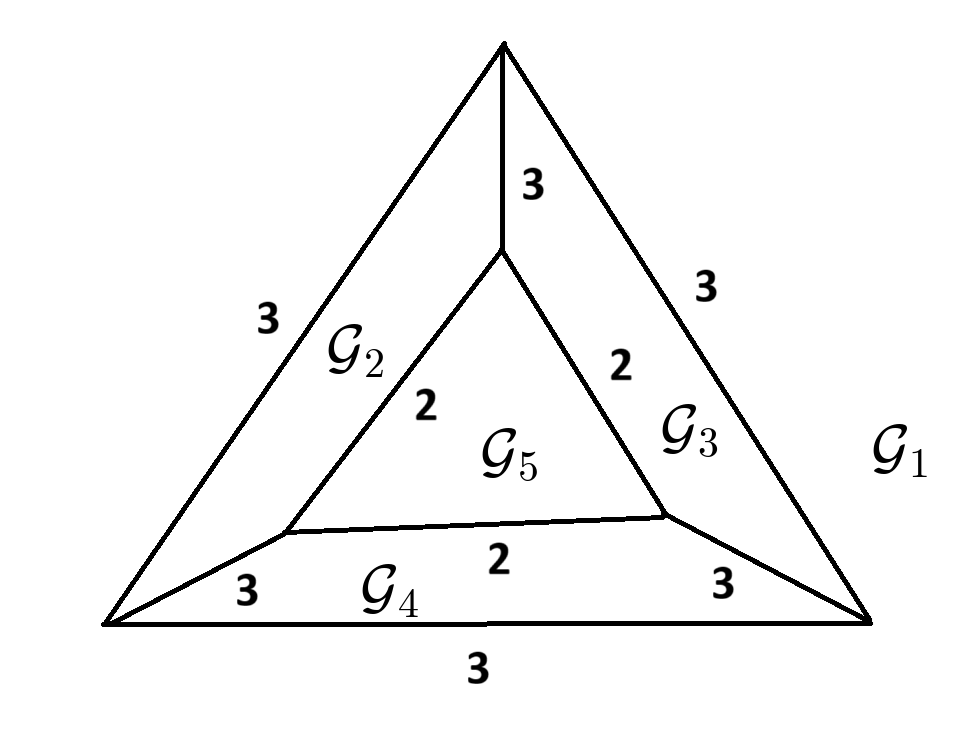}
		\caption{The truncation polytope used in Example~\ref{ex:prism-one-vertex-truncation}.}
		\label{fig:prism}
	\end{figure}

Let \(\mathcal G\) be the triangular prism obtained from a tetrahedron by truncating one vertex.
We denote the truncation facet by \(\mathcal G_5\) and the opposite triangular facet by \(\mathcal G_1\).
The remaining facets \(\mathcal G_2,\mathcal G_3,\mathcal G_4\) are quadrilaterals.
The three edges of \(\mathcal G_5\) are labeled by \(2\), so \(m_{25}=m_{35}=m_{45}=2\).
We label all remaining edges by \(3\).
The only nonadjacent pair is \((\mathcal G_1,\mathcal G_5)\), and hence \(m_{15}=\infty\).
See Figure~\ref{fig:prism}.

Since \(\Sigma(\mathcal G)\) is connected, Remark~\ref{rem:truncation-group-hypotheses} shows that \(\pi_{\mathcal G}\) is infinite, irreducible, and non-affine.
Let \(C=(c_{ij})\) be a Cartan matrix of a \(\mathcal G\)-compatible representation into \(\SLpm{5}\).
After normalizing \(C\) along the directed spanning tree \(5\to1,\;1\to2,\;2\to3,\;3\to4\), with all prescribed entries equal to \(-1\), we obtain
\[
C=
\begin{pmatrix}
	2&-1&-1/z&-1/w&-(x+4)\\
	-1&2&-1&-y&0\\
	-z&-1&2&-1&0\\
	-w&-1/y&-1&2&0\\
	-1&0&0&0&2
\end{pmatrix}.
\]
Here \(x\ge0\) and \(y,z,w>0\).
These parameters are uniquely determined by the normalized matrix.
Conversely, every matrix of this form is compatible with \(M(\mathcal G)\) and has rank at most five.
By Danciger et al.~\cite[Facts~3.28--3.29]{danciger2025convex} and Lemma~\ref{lem:automatic-v5}, the continuous normalization above therefore identifies \(\Xref{5}(\mathcal G)\) homeomorphically with the domain \(x\ge0,\ y,z,w>0\).
In this normalization, the condition \(\mu_{15}>4\) is equivalent to \(x>0\).
Hence \(\Xstr{5}(\mathcal G)\) is homeomorphic to the domain \(x,y,z,w>0\).

Define
\[
\begin{aligned}
	N(y,z,w)={}&
	w+2wy+2w^2y+w^2y^2+2z+yz+9wyz+w^2yz\\
	&+2w^2y^2z+z^2+2yz^2+2wyz^2+wy^2z^2.
\end{aligned}
\]
A direct computation gives
\[
\det C=\frac{wz(1-y)^2x-2N(y,z,w)}{wyz}.
\]
If \(y=1\), then \(\det C<0\).
For \(y\ne1\), the equation \(\det C=0\) is equivalent to
\[
x=\frac{2N(y,z,w)}{w(1-y)^2z}.
\]
Since \(N(y,z,w)>0\) for \(y,z,w>0\), this equation defines a positive continuous function \(x\) on each of the domains \(\{0<y<1,\ z,w>0\}\) and \(\{y>1,\ z,w>0\}\).

By Lemma~\ref{lem:strict-character-rank-bounds}, a character in \(\Xstr{5}(\mathcal G)\) has a Cartan matrix of rank four if and only if its normalized Cartan matrix satisfies \(\det C=0\).
Under the parametrization by Cartan matrices, the two graphs defined by \(\det C=0\) are the connected components of \(X^{\mathrm{ref},(4)}_{5,\mathrm{str}}(\mathcal G)\), each homeomorphic to \(\mathbb R^3\).

The rank-five subspace corresponds to the disjoint union of the following three open regions in the strict parameter domain:
\[
\{\det C<0\},\qquad
\{\det C>0,\ 0<y<1\},\qquad
\{\det C>0,\ y>1\}.
\]
The two positive-determinant regions are connected, since each is the region above one of the two continuous graphs.
The negative-determinant region is connected, since along any compact line segment in the positive \((y,z,w)\)-domain, the determinant is negative for all sufficiently small \(x>0\).
These are therefore the three connected components of \(X^{\mathrm{ref},(5)}_{5,\mathrm{str}}(\mathcal G)\).
\end{example}

\begin{example}\label{ex:non-coxeter-realizable-empty-strict}
	\begin{figure}[htbp]
		\centering
		\includegraphics[width=0.36\textwidth]{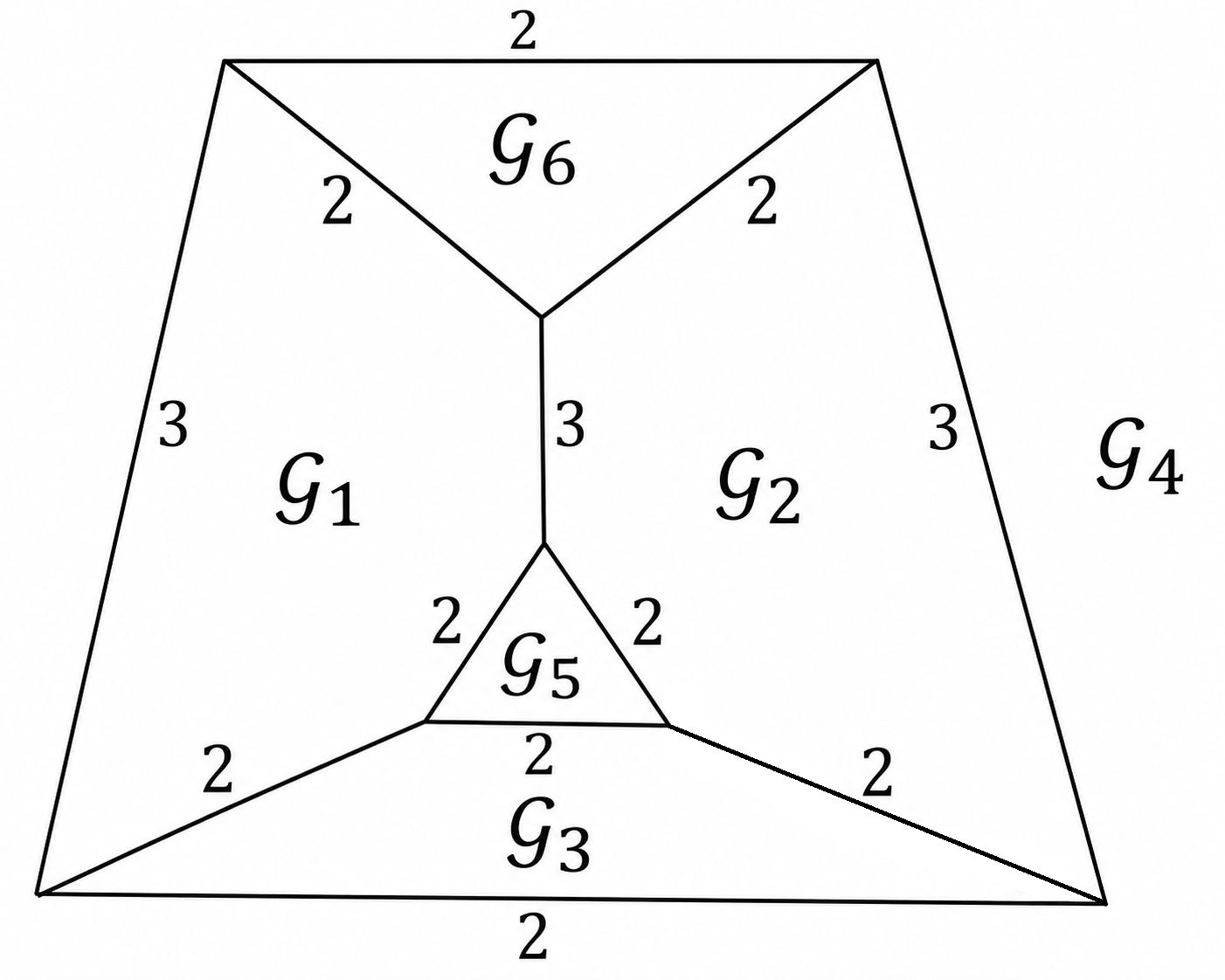}
		\caption{The truncation polytope used in Example~\ref{ex:non-coxeter-realizable-empty-strict}.}
		\label{fig:two-adjacent-vertex-truncations}
	\end{figure}

	This example exhibits a labeled truncation \(3\)-polytope whose full reflection character variety in \(\SLpm{5}\) is nonempty and whose subspace of strict characters is empty.
	Let \(\mathcal G\) be obtained from a tetrahedron by truncating two adjacent vertices.
	Label the four facets inherited from the tetrahedron by \(\mathcal G_1,\mathcal G_2,\mathcal G_3,\mathcal G_4\).
	Choose the marking so that the truncated vertices are opposite \(\mathcal G_4\) and \(\mathcal G_3\) in the original tetrahedron.
	Denote the corresponding truncation facets by \(\mathcal G_5\) and \(\mathcal G_6\), respectively.
	Thus \(\mathcal G_5\) is adjacent to \(\mathcal G_1,\mathcal G_2,\mathcal G_3\), and \(\mathcal G_6\) is adjacent to \(\mathcal G_1,\mathcal G_2,\mathcal G_4\).
	The nonadjacent pairs are \((\mathcal G_4,\mathcal G_5)\), \((\mathcal G_3,\mathcal G_6)\), and \((\mathcal G_5,\mathcal G_6)\).
	We assign Coxeter labels
	\[
	m_{15}=m_{25}=m_{35}=m_{16}=m_{26}=m_{46}=2,
	\]
	and
	\[
	m_{12}=m_{14}=m_{24}=3,\quad m_{13}=m_{23}=m_{34}=2,
	\]
	as in Figure~\ref{fig:two-adjacent-vertex-truncations}.

	The resulting \(\mathcal G\) is not Coxeter realizable in \(\mathbb S^3\).
	Indeed, \(\mathcal G_1\), \(\mathcal G_2\), and \(\mathcal G_3\) are pairwise adjacent but have no common vertex, since the truncation creating \(\mathcal G_5\) removes their common intersection.
	The prescribed dihedral angles between these facets are \(\pi/3\), \(\pi/2\), and \(\pi/2\), whose sum exceeds \(\pi\).
	Marquis~\cite[Theorem~3.16\textup{(3b)}]{marquis2010espace} therefore gives \(\Ccal(\mathcal G)=\varnothing\).
	
	We next show that \(\Xref{5}(\mathcal G)\) is nonempty while \(\Xstr{5}(\mathcal G)\) is empty.
	Since \(\Sigma(\mathcal G)\) is connected, Remark~\ref{rem:truncation-group-hypotheses} shows that \(\pi_{\mathcal G}\) is infinite, irreducible, and non-affine.
	For every diagonal equivalence class of Cartan matrices compatible with \(M(\mathcal G)\), normalization along the directed spanning tree \(1\to2,\;2\to4,\;4\to5,\;5\to6,\;6\to3\), with all prescribed entries equal to \(-1\), gives a representative of the form
	\[
	C=
	\begin{pmatrix}
		2&-1&0&-x&0&0\\
		-1&2&0&-1&0&0\\
		0&0&2&0&0&-z\\
		-1/x&-1&0&2&-1&0\\
		0&0&0&-y&2&-1\\
		0&0&-1&0&-w&2
	\end{pmatrix},
	\]
	where \(x>0\), \(y=\mu_{45}\ge4\), \(z=\mu_{36}\ge4\), and \(w=\mu_{56}\ge4\).
	The only edge of \(\Sigma(\mathcal G)\) outside this tree is \(14\), with entries \(c_{14}=-x\) and \(c_{41}=-1/x\).

	If \(C\) arises from a reflection pair in \(\mathbb R^5\), then \(\operatorname{rank}C\le5\), and hence \(\det C=0\).
	A direct computation gives
	\[
	x\det C=(z-4)\bigl(3xy+2(x-1)^2\bigr)+2w(x-1)^2.
	\]
	Since \(3xy+2(x-1)^2>0\), \(z\ge4\), and \(w>0\), we have \(\det C=0\) if and only if \(z=4\) and \mbox{\(x=1\)}.
	Consequently, every representation in \(\Homrefc{\mathcal G}(\pi_{\mathcal G},\SLpm{5})\) satisfies \(\mu_{36}=4\).
	It follows that \(\Homstrictc{\mathcal G}(\pi_{\mathcal G},\SLpm{5})=\varnothing\), and hence \(\Xstr{5}(\mathcal G)=\varnothing\).

	The full reflection character variety is nevertheless nonempty.
	For example, take \(y=w=5\), \(z=4\), and \(x=1\).
	The displayed matrix is compatible with \(M(\mathcal G)\) and has rank at most five.
	Danciger et al.~\cite[Fact~3.28]{danciger2025convex} therefore give a semisimple representation generated by weakly compatible reflections on \(\mathbb R^5\) with this Cartan matrix.
	Lemma~\ref{lem:automatic-v5} supplies \textup{(V5)}, so this representation is \(\mathcal G\)-compatible.
	Hence \(\Xref{5}(\mathcal G)\ne\varnothing\) while \(\Xstr{5}(\mathcal G)=\varnothing\), so the nonemptiness hypothesis in Theorem~\ref{thm:main-dimension}\textup{(iii)} is essential.
\end{example}

\begin{example}\label{ex:non-two-perfect-dimension}
	This example shows that the \(2\)-perfectness assumption in Theorem~\ref{thm:main-dimension}\textup{(iii)} is essential.
	Let \(\mathcal G\) be obtained from a \(5\)-simplex by truncating one vertex.
	Denote the truncation facet by \(\mathcal G_7\) and the facet opposite the truncated vertex in the original simplex by \(\mathcal G_1\).
	Mark the remaining facets by \(\mathcal G_2,\ldots,\mathcal G_6\).
	We assign the labels
	\[
	m_{ij}=2\qquad(i\in\{1,2,3\},\ j\in\{4,5,6\})
	\]
	and label all remaining ridges by \(3\).
	The only nonadjacent pair is \((\mathcal G_1,\mathcal G_7)\), so \(m_{17}=\infty\).
	Thus \(d=5\), \(f=7\), and \(e_2=9\).
	The facets \(\mathcal G_1,\ldots,\mathcal G_4\) meet in an edge \(E\).
	The Coxeter diagram of \(W_{I_E}\) contains a triangle with Coxeter labels \(3,3,3\), so \(W_{I_E}\) is infinite.
	Hence \(\mathcal G\) is not \(2\)-perfect.
	
	Since \(\Sigma(\mathcal G)\) is connected, Remark~\ref{rem:truncation-group-hypotheses} shows that \(\pi_{\mathcal G}\) is infinite, irreducible, and non-affine.
	Let \(C=(c_{ij})\) be a Cartan matrix of a \(\mathcal G\)-compatible representation into \(\SLpm{6}\).
	After normalizing \(C\) along the directed spanning tree \(7\to1,\;1\to2,\;2\to3,\;7\to6,\;6\to5,\;5\to4\), with all prescribed entries equal to \(-1\), we obtain
	\[
	C=
	\begin{pmatrix}
		2&-1&-x&0&0&0&-z\\
		-1&2&-1&0&0&0&-p\\
		-1/x&-1&2&0&0&0&-q\\
		0&0&0&2&-1&-y&-r\\
		0&0&0&-1&2&-1&-s\\
		0&0&0&-1/y&-1&2&-1\\
		-1&-1/p&-1/q&-1/r&-1/s&-1&2
	\end{pmatrix},
	\]
	where \(x,y,p,q,r,s>0\) and \(z=\mu_{17}\ge4\).
	The principal submatrices of \(C\) on \(\{1,2,3\}\) and \(\{4,5,6\}\) have determinants \(2-x-x^{-1}\le0\) and \(2-y-y^{-1}\le0\), respectively.
	Their adjugates have positive entries.
	A block determinant expansion therefore gives
	\[
	\det C=0 \;\Longleftrightarrow \; x=y=1.
	\]
	Since \(\rank C\le6\), we obtain \(x=y=1\).
	For these values, solving \(C\xi=0\) shows that \(\ker C\) is one-dimensional, so \(\rank C=6\).
	
	Conversely, every matrix of this form with \(x=y=1\), \(p,q,r,s>0\), and \(z\ge4\) is compatible with \(M(\mathcal G)\) and has rank six.
	By Danciger et al.~\cite[Facts~3.28--3.29]{danciger2025convex} and Lemma~\ref{lem:automatic-v5}, the continuous normalization above therefore identifies \(\Xref{6}(\mathcal G)\) homeomorphically with the domain
	\[
	\{(p,q,r,s,z)\in\mathbb R_{>0}^5:z\ge4\}.
	\]
	Under this identification, strict characters correspond to \(z>4\).
	Consequently,
	\[
	\Xstr{6}(\mathcal G)\ne\varnothing,
	\qquad
	\dim\Xref{6}(\mathcal G)=\dim\Xstr{6}(\mathcal G)=5.
	\]
	However, the dimension formula in Theorem~\ref{thm:main-dimension}\textup{(iii)} gives
	\[
	(2\cdot6-5-2)\cdot7-6^2+\frac{5\cdot6}{2}+1-9=6.
	\]
	Thus the dimension equality fails without the \(2\)-perfectness assumption.
\end{example}

	\begingroup
	\setlength{\emergencystretch}{1.25em}
	\printbibliography[heading=bibintoc, title={References}]
	\endgroup
\end{document}